\documentclass[a4paper,11pt]{article}

\RequirePackage[OT1]{fontenc}
\RequirePackage{amsthm,amsmath}
\RequirePackage[numbers]{natbib}
\RequirePackage[colorlinks,citecolor=blue,urlcolor=blue]{hyperref}
\RequirePackage{hypernat}

\usepackage[applemac]{inputenc}
\usepackage[french]{babel}

\usepackage{latexsym,enumerate}
\usepackage{amsmath,amsthm,amsopn,amstext,amscd,amsfonts,amssymb,amsbsy}
\usepackage{mathrsfs}
\usepackage{fullpage}
\usepackage{graphicx}
\usepackage[usenames]{color}

\newtheorem{proposition}{Proposition}
\newtheorem{lemma}{Lemma}
\newtheorem{theorem}{Theorem}

\newtheorem{corollary}{Corollary}

\let\epsilon=\varepsilon
\let\eps=\epsilon
\let\phi=\varphi

\let\tilde=\widetilde

\newcommand{\field}[1]{\mathbb{#1}}
\newcommand{\R}{\field{R}}

\newcommand{\cI}{\mathcal{I}}
\newcommand{\cL}{\mathcal{L}}

\newcommand{\cJ}{{\mathcal J}}

\newcommand{\beqn}{\begin{equation}}
\newcommand{\eeqn}{\end{equation}}
\newcommand\eref[1]{(\ref{#1})}

\def\I{{\field I}}

\begin{document}
\begin{center}
{\sc \Large A new selection method for high-dimensionial
instrumental setting: application to the Growth Rate convergence
hypothesis} \vspace{0.2cm}

\vspace{1cm}

Mathilde MOUGEOT $^{\mbox{\footnotesize a}}$, Dominique
PICARD$^{\mbox{\footnotesize a}}$ and Karine
TRIBOULEY$^{\mbox{\footnotesize b}}$

\vspace{0.5cm}

$^{\mbox{\footnotesize a}}$ LPMA\\Universit\'e Diderot -- Paris VII\\   175 rue du Chevaleret\\
75013 Paris, France\\
\smallskip
\textsf{mathilde.mougeot@univ-paris-diderot.fr}\\
\textsf{picard@math.jussieu.fr}\\
\bigskip
$^{\mbox{\footnotesize b}}$ Modal'X, Universit\'e Paris Ouest --  Paris X,\\
200 rue de la R\'epublique, 92000 Nanterre, France\\
\smallskip
\textsf{karine.tribouley@u-paris10.fr}\\

\vspace{0.5cm}

\end{center}

\vspace{0.5cm}

\begin{abstract}
\noindent {\rm This paper investigates the problem of selecting
variables in regression-type models for an "instrumental" setting.
Our study is motivated by empirically verifying the conditional
convergence hypothesis used in the economical literature concerning
the growth rate. To avoid unnecessary discussion about the choice
and the pertinence of instrumental variables, we embed the model in
a very high dimensional setting. We propose a selection procedure
with no optimization step called LOLA, for Learning Out of Leaders
with Adaptation. LOLA is an auto-driven algorithm with two
thresholding steps. The consistency of the procedure is proved under
sparsity conditions and simulations are conducted to illustrate the
practical good performances of LOLA.  The behavior of the algorithm
is studied when instrumental variables are artificially added
without a priori significant connection to the model. Using our
algorithm, we provide a solution for modeling the link between the
growth rate and the initial level of the gross domestic product and
empirically prove the convergence
hypothesis. \\

}
\end{abstract}

\newpage

\baselineskip=18 pt

\section{Introduction}
Let us consider the usual linear model $Y=X\alpha+u$ stated to
explain the causal relationship of explanatory variables $X$ on
another variable $Y$ in the specific and problematic case where the
p-value associated to the Fisher test is small implying that the
regression is not significant. This can happen if the covariates $X$
are endogenous: for  instance when determinant predictors are
omitted from the model or when the covariates $X$ are correlated
with the errors $u$. The insertion of instrumental variables $Z$  in
the model may lead to
 consistent estimation of the coefficients $\alpha$. In
economy theory, variables reporting behavior breaks (policy and
political changes, natural disasters) appear to be good candidates
to create instrumental variables. Choosing the instruments is a
serious issue. We propose to address this issue in an objective way
by considering a huge set of potential variables $Z_1,\ldots,Z_p$.
Then obviously, a new problem appears: selecting among this huge set
of potential candidates the relevant ones. When the number $p$ of
candidates is very large with respect to the number $n$ of available
observations, the usual classical selection methods fail. The goal
of our paper is to present a very simple procedure called LOLA (with
no recursive step and no optimization step) dedicated to high
dimensional regression models and to examine its selection
properties when confronted to this specific problem of selecting
instrumental variables.

Let us first present the specific question concerning the
International Economic Growth. One especially important point in the
empirical growth literature is to evaluate the effect of an initial
level $X$ of gross domestic product (GDP) per capita on the growth
rate of GDP $Y$. The conditional convergence hypothesis states that,
other things being equaled, countries with lower GDP per capita are
expected to grow more than others due to higher marginal returns on
capital stock. This convergence hypothesis should imply a negative
effect. However, when empirical tests are performed on the simple
linear model $Y=\alpha_1+\alpha_2X+u$, the alternative "$H_1:\;
\alpha_2<0$" is rejected (see Table 2, Section~\ref{PNBlola}). A
natural idea is that other phenomena interfere in this relationship
hiding the negative effect. Unfortunately, growth theory is not
 explicit enough about the set of key factors of the growth. These
factors could belong to various categories: policy variables
(fiscal, exchange rate and trade policies), political variables
(rule of law, political rights, ...), religious variables, regional
variables, type of investment (equipment/non equipment), variables
relating to the macroeconomic environment (inflation, initial level
of GDP,...), variables accounting for the international environment
(terms of trade, ...), etc. Let us explain for instance how Sala
(1997) proceeds to select the instrument variables among all the
potential factors. Based on previous studies, three predictors (the
level of income, the life expectancy and the primary-school
enrollment rate) are retained a priori.  Next, all the possibilities
of introducing
 three other predictors (the number of variables is then
$p=1+3+3=7$) are inspected  and
the regression model is estimated.
 A test is then built to choose among  all
the estimated regression models. The methodology explains the title
of the paper: "I just ran two millions regression". We present here
a different approach: we consider  a huge  quantity of predictors to
be selected and we compute only ONE regression model but in high
dimension.

The growth rate problem is also studied in Belloni and Chernozhukov
(2010) where results are obtained using the lasso methodologies
implying optimization steps. Since their results are convincing, we
think the methodology is promising and we compare our results to
theirs. We also add a nonparametric point of view, which shed a new
light on the construction of instrumental variables. Instead of
adding economic variables, we just consider the data as a signal and
analyze it as depending on two factors: the initial level $X$ of
gross domestic product plus a unknown signal function estimated in a
nonparametric way. As in Belloni and Chernozhukov (2010), we use the
Barro and Lee data base which contains a huge number of instruments
$Z=(Z_1,\ldots,Z_p)$ including many covariates for characterizing
the different countries (see Section~\ref{PNB} and for a detailed
description of the data, see Barro and Sala (1995) ). Since the
number of instruments is large $p\sim 300$ compared to the number of
countries $n\sim 100$, we are in a high dimensional setting and we
proceed as follows:

$\bullet$ First, we apply the selection procedure LOLA  to
reduce the dimensionality of the problem and to extract the relevant
instruments. This procedure is auto-driven which has the advantage
to avoid the search for tuning parameters. Let us emphasize that the number
${\widehat{S}}$ of selected instruments is an output of our
procedure and is not imposed a priori.

$\bullet$  Next, we use the selected instruments to estimate the
 parameter $\alpha_2$ describing the relationship between $X$ (GDP) and $Y$ (growth rate).

The results obtained with LOLA are excellent to verify empirically
the conditional convergence hypothesis: the estimation of parameter
$\alpha_2$ is negative and roughly speaking, the null hypothesis
"$H_0:\,\alpha_2=0$" is rejected even for small significant
prescribed test levels. Some determinants for the growth rate are
identified and we notice differences
 according to the periods of time considered in the study.
Moreover, some selected variables are the same as the instruments
selected in Belloni and Chernozhukov (2010) even though we consider
a larger number of predictors.

 \vspace{0.6cm}

Let us come back now to the more general setting $Y=X\beta+u$ where
$X=(X_1,\ldots,X_p)$ is a list of $p$ predicting variables. These
models are situated in a stream of papers devoted to the estimation
of $\beta$ as well as the selection in a high dimensional setting.
The most accomplished methods in this context have in common a
crucial assumption called "Sparsity Assumption". Roughly speaking,
it is assumed that, even if the model depends on a huge number of
parameters, only a (very) small number of them are significant.
Hence a selection step is conceivable. Basically, the numerous
methods which are proposed in this context can be classified in
three categories: the filtering methods retaining the best
explanatory variables (using various criteria), the 'wrapping'
methods (among them the greedy algorithms) operating a selection in
a stepwise way, and the optimization methods (such as Dantzig
selector or Lasso minimisation) minimizing a criteria combining an
empirical risk with a penalization term related to the number of
retained coefficients. Although it is impossible to be exhaustive in
such a productive domain, we cite among many others Fan and Lv
(2010, 2008), Tibshirani (1996), Candes and Tao (2007), Bunea,
Tsybakov,
 and Wegkamp (2007), Needell and Tropp (2009), Tropp and Guilbert (2007),
Barron, Cohen,  Dahmen and DeVore (2008), Haury, {Jacob}, and {Vert}
(2010) (and apologize for all the papers which would deserve to be
cited here).
 Since it is a rough common feeling in the applied community that these methods are generally more taylored to  prediction, often generate instability in the selection step, and that the simplest ones (filtering) finally
 behave very honorably for selection criteria, we detail here a selection
 method (LOL), which situates in between filtering and wrapping methods,
 since it has only two selection steps, but remains one of the simplest one.
 This methods has been proved to have, despite its simplicity, theoretical
 optimal properties for the prediction criteria
(see Kerkyacharian,Mougeot,Picard and Tribouley (2009) and Mougeot,
Picard and Tribouley (2010) ).

 In this paper, we investigate the selection properties
  of this procedure. From a theoretical point
  of view, we state exponential convergence of the
  false positive and false negative rates under fairly general conditions.
  Then, our main focus is to study the practical performances of the selection procedure
  in the particular domain of instrumental variables.
Let us
describe in more details LOL selection procedure. Two consecutive
steps of thresholding are performed. In each of these steps we
'kill' variables and the result of our selection is the variables
which have successfully passed the two steps.

$\bullet$   The first step is a thresholding procedure allowing to
reduce the dimensionality of the problem in a rather rough way by a
simple inspection of the "correlations", computed between the target
variable and the predictors.

$\bullet$  The least square method is then used on the linear
sub-model defined by  the variables retained after the first step.
The second step of thresholding is performed on the estimations of
the parameters of the sub-model. This  step is more refined and
corresponds to a denoising phase of the algorithm.

\vspace{0.3cm}

This procedure is called LOL for "Learning Out of Leaders". The
thresholding levels $t$ and $s$ used in  both steps are the inputs
of  LOL algorithm and are set by the statistician. Theoretical
results are established in terms of $s$ and $t$ and more precisely,
Proposition \ref{false_detection} states the consistency of the LOL
algorithm in the sense that the number of false detections as well
as the number of false negative tend to zero when the number of
observations tends to infinity. Some assumptions are obviously
needed to obtain the convergence properties: the {\textit{sparsity
assumption}} (only a few parameters are significant even if we do
not know their number and position), the {\textit{significance
assumption}} (the 'significant' coefficients are above a
'significant' level). Both properties are standard, particularly in
the high dimensional setting. We address also in an experimental way
some problems which are not solved in the theoretical part. Among
them, we focus on building an operational procedure which is auto
driven (i.e.  the levels $s$ and $t$ are chosen in an adaptive way).
This procedure is called
 LOLA procedure for LOL completed by an A algorithm allowing  adaptation.
We provide an illustration of LOLA thanks to simulations considering
the very common case where the predictors are Gaussian variables.
We first observe that LOLA appears to be an extremely accurate procedure when the
predictors are independent or when the number of predictors to be
selected is small. Hence, in a second step, we relax  these ideal conditions by
considering dependent variables and  increasing the number of
variables to be selected. We observe that the results given by LOLA
procedure are still very convincing.
We also investigate the performance of LOLA in a toy instrumental setting related to the Boston Housing data. Finally, we address, using Baro
and Lee data, the convergence hypothesis problem in economy.

The paper is organized as follows. In Section~\ref{section_algo},
LOL algorithm is described , the general model and the hypotheses on
the model are presented. Theoretical results on the
consistency of the LOL procedure is
 stated in Theorem \ref{consistence_MS}. In Section~\ref{section adaption}, we explain
how to modify LOL into LOLA to obtain a data driven procedure and we
explore practical performances of LOLA with some simulations.
Section~\ref{Boston} is dedicated to explore a first toy example
using a classical dataset (Boston Housing) with a huge additional
set of simulated 'instrumental' variables. Before applying to real
data, the purpose of this section is to verify the ability of our
algorithm to accurately select the appropriate variables, when an
important set of unappropriate variables are added. Finally, in
Section~\ref{PNB} we focus on the central question and we prove the
hypothesis convergence used in the Solow-Swan-Ramsey growth model.
The proofs of the theoretical results are detailed in
Section~\ref{section_proof}.

\section{LOL and  theoretical properties}\label{section_algo}
The selection procedure need two tuning parameters which are inputs
given by the user. The assumptions needed to obtain the theoretical
results are presented in this case. The consistency of the procedure
is established in Theorem \ref{consistence_MS}.
\subsection{LOL Procedure}\label{essentialLOL}
The selection algorithm is denoted LOL$(X,Y,t,s)$ : the inputs
 are the target variable $Y$, the predictor variables
$(X_1,\ldots,X_p)$ and  the two tuning parameters $t,s$, specified by the user. As output, the procedure
provides the set  of indices $\widehat{\cI}$ of the selected variables.

\vspace{0.5cm}

\begin{center}
\label{algolol}
\begin{tabular}{|l|ll|}\hline
{\bf $\widehat{\cI}\leftarrow$ LOL$(X,Y,t,s)$}& {\bf Input:}
&target  $Y$, regression
variables $X=(X_1,\ldots,X_p)$\\
&&
tuning parameters $t,s$\\
&{\bf Output:} &set of indices of the selected regression variables
$\widehat{\cI}$\\
 \hline\end{tabular}\end{center}

 \vspace{0.5cm}

Let us describe LOL$(X,Y,t,s)$. The regression variables $X$ are first normalized
and defined by $\tilde{X}=XD$ where $D$ is a $(p\times p)$
diagonal matrix $D_{\ell m}:= \sigma_\ell^{-1} \,\I_{\ell=m}$ and
 $\sigma^2_\ell=n^{-1}(X_\ell^t X_\ell)$ is the empirical
second moment of the predictor $X_\ell$. The coherence of the matrix
of the normalized predictors is then computed:
$$\tau_n=\sup_{\ell\not=m}|\frac1n\sum_{i=1}^n
\tilde{X_{i\ell}}\tilde{X_{im}}|$$ The 'correlations' (scalar products) between the
target variable $Y$ and the normalized predictors $\tilde{X_\ell}'$s
are then sorted by increasing order:
$$
|\tilde{X_{\bullet}}^tY|_{(1)}\geq
|\tilde{X_{\bullet}}^tY|_{(2)}\geq \ldots \geq
|\tilde{X_{\bullet}}^tY|_{(p)}.
$$
The leaders are the predictors associated with the highest
'correlations', with indices belonging to the set
\begin{eqnarray}\label{leader}
\cJ=\,\left\{ \ell=1,\ldots,p,\;
 |\tilde{X_{\ell}}^tY|\geq \left(|\tilde{X_{\bullet}}^tY|_{
 (\lfloor\tau_n^{-1}\rfloor)}\vee t\right)\right\}
\end{eqnarray}
where $t$ is one of the input of LOL algorithm . Denoting the
extracted matrix $\tilde{X_{\cJ}}$ by
\begin{eqnarray*}(\tilde{X_\cJ})_{i,\ell} =
 \tilde{X_{i\ell}}\quad \mbox{ for any } \ell\in\cJ \mbox{ and }
 i\in\{1,\ldots,n\}
 \end{eqnarray*}
the Ordinary Least Square (OLS) estimator for the pseudo linear model
$$\forall \ell=1,\ldots,p,\quad \widehat{ \beta_\ell}=
\left((\tilde{X_\cJ}^t \tilde{X_\cJ})^{-1}\tilde{X_\cJ}^t
\,Y\right)_\ell\,\I\{ \ell\;\in \cJ\}$$ the set of indices of the
selected predictors is finally given by
$$
\widehat{\cI}=\{\ell=1,\ldots,p,\; |\widehat{\beta_\ell}|\ge s\}
$$
where $s$ is the second parameter input of the  LOL algorithm.

\subsection{Model assumptions and selection criteria\label{model}}
We consider a Gaussian (or sub-gaussian) high dimensional linear
model. More precisely, we assume that the target variable $Y$ and the
$p$ predictors $X=(X_1,\ldots,X_p)$ are
 linked through the linear regression model
\begin{align*}
 Y&=X\beta+\eps
 \end{align*}
 where  $\beta\in \R^p$ is an unknown parameter
   and
the vector $\eps= (\eps_i)_{i=1,\ldots,n} $
  is a vector of  independent Gaussian variables  $N(0,\eta^2)$.
In the selection problem, it is generally assumed
 that only a few coefficients $\beta_\ell$ are non zero. The set of
 non zero coefficients
 $$\cI:=\left\{\ell\in\; \{1,\ldots,p\},\beta_\ell\not=0\right\}$$
 is precisely
 the set to be estimated.
A sparsity condition is introduced to enforce the cardinality of $\cI$
to be less than $S$. In practice such kind of assumption is not
realistic so  a relaxed setting is here considered:  only a few coefficients are
  larger than a first threshold $s_n$
  and the significant coefficients (the ones we definitely want to detect)
  are larger than a second threshold $t_n$.
More precisely, we assume:
 \begin{itemize}
\item {\bf Sparsity Conditions:} There exist $S>0$ and thresholds
levels $t_n$ and $s_n$ such that
\begin{align}
 \#\left\{\ell\in\{1,\ldots,p\},\;|\sigma_\ell\beta_\ell|\ge s_n/2\right\}\;\le{S}
 \quad \mbox{ and }\quad
 \sum_{\ell=1}^p|\sigma_\ell\beta_\ell|^2
 \;I\{|\sigma_\ell\beta_\ell |\le 2t_n\}\le \,\frac {S\log{p}}n
 \label{Scond}
 \end{align}
$\sigma_l$ are the normalizing factors defined in Section
\ref{essentialLOL}.
 \item {\bf Size Condition:} There exist a positive constant $M$ such that \begin{align}
\sum_{\ell=1}^p|\sigma_\ell\beta_\ell|\le { M} \label{B1cond}
 \end{align}
 \item{\bf Significance:}
 There exists   a sequence $\mu_n>0$ such that
 the set $\cI$ of coefficients to be detected is defined by
\begin{align}
 \cI:= \left\{\ell\in\{1,\ldots,p\},\; |\sigma_\ell\beta_\ell|\geq \mu_n\right\}\label{B2cond}
 \end{align}
  \end{itemize}
  Recall that a selection procedure
whose output is denoted by $\widetilde{\cI}$ is said to be
consistent if $P(\,\widetilde{\cI}=\cI\,)$ is tending to $1$ when
the number $n$ of observations is tending to infinity. More
precisely, to evaluate the quality $\widetilde{\cI}$, the number  of
False Positive $FP$ and  the number of False Negative decisions
$FN$ are generally computed:
$$FP:=\sum_{\ell =1}^pI\{\ell \not\in\cI\}I\{\ell
\in \widetilde{\cI}\}\quad \mbox{ and }\quad FN:=
\sum_{\ell=1}^pI\{\ell\in\cI\}I\{\ell\not \in \widetilde{\cI}\} .$$

\subsection{Performances}\label{result}
The performances of LOL algorithm are depending - and this is a
common feature in the regression problem - on  the regression matrix
$X$ and particularly on the coherence  $\tau_n$
 of the $(p\times
p)$ matrix $(D^tXXD)$ defined as
\begin{equation*}
\tau_n=\sup_{\ell\not=m }|\,(DX^tXD)_{\ell m}|
\end{equation*}
where $D$ is the $(p\times p)$ diagonal matrix $D_{\ell m}:=
\sigma_\ell^{-1} \,\I_{\ell=m}$ and
$\sigma_\ell^2:=\frac1n\sum_{i=1}^n X_{i\ell}^2$ is the empirical
second moment of variable $X_\ell$.
 This quantity is important because it induces a
bound on the size of
 the invertible sub-matrices built with the columns of $(D^tXXD)$ and
 thus on the maximum number of leaders used by our algorithm.
 For more details, we refer to Mougeot, Picard and Tribouley (2010).

The consistency results are first stated for general threshold
levels $t_n,s_n$.
 \begin{theorem}\label{consistence_MS}
Suppose that $p\le \exp (an)\; $ for some constant $a>0$.
 Choose the thresholds $t_n, s_n$ such that $t_n\geq s_n\vee\tau_n\vee
 \sqrt{\log{p}/n}$. Recall that $\widehat \cI$ is the output of the algorithm LOL($X,Y,t_n,s_n$). Then
\begin{align*}
P\left( \cI  = \widehat{\cI}\right) \geq 1-
\exp\left(-c\,ns_n^2\right)
\end{align*}
as soon as the vector $\beta$ verifies the sparsity conditions
(\ref{Scond}) as well as the size and significance conditions (\ref{B1cond}), (\ref{B2cond}) for
  $$\mu_n\geq O\left(t_n \vee\tau_ns_n\,\sqrt{n} \vee
  s_n\,\sqrt{|\log{\tau_n}|}\right)\quad\mbox{
  and }\quad S\leq O(ns_n^2).$$
\end{theorem}

Theorem \ref{consistence_MS} is a consequence of the following
proposition which gives distinct evaluations of the errors induced
by the false positive detections and  the false negative detections.
It is interesting to observe that slightly different assumptions are
needed for both detections: for instance, no explicit condition for
$\mu_n$ is needed to insure the convergence of the rate of False
Negative detections.

\begin{proposition}\label{false_detection} Let $k$ be a given positive
number. Assume that $p\le \exp( an)\; $ for some $a>0$.
\begin{itemize}
\item {\bf False Negative (FN):} Choose the threshold levels $t_n\geq s_n$
and assume that the vector $\beta$  verifies the sparsity conditions
(\ref{Scond}) as well as the size and significance conditions (\ref{B1cond}), (\ref{B2cond}) for

  \begin{align}\mu_n\geq O\left(t_n \vee\tau_n\sqrt{\frac{ S}{k}} \vee
\sqrt{\frac{S}{nk}}\vee \sqrt{\frac{S|\log{\tau_n}|}{nk}}\right).
\label{num1}
\end{align}
Then there exists a constant $c>0$ such that
\begin{align*}
P(FN>k)&\leq \exp\left(-c\,kn\mu^2_n\right).
\end{align*}
\item {\bf False Positive (FP):}
Choose the thresholds
\begin{align}t_n\geq s_n\vee O\left( \sqrt{\log{p}/n}\vee
\,\tau_n \right)\label{num3}
\end{align}
and assume that the vector $\beta$  verifies the sparsity conditions
(\ref{Scond}) as well as the size and significance conditions (\ref{B1cond}), (\ref{B2cond}) for
 $S\leq
O(nks_n^2)$.

 Then there exists some constant $c>0$ such that
\begin{align}
P(FP>k)&\leq \exp\left(-c\,kn s^2_n\right).\label{num2}
\end{align}
\end{itemize}
\end{proposition}

Observe that the requirements on the thresholds are quite loose but
clearly the performances of the algorithm suffer when the thresholds
are too low (by weakening the exponential rates -see \eref{num2}, \eref{num3}-). They also suffer
when the thresholds are too high since the significant coefficients
have to be above the thresholds -see \eref{num1}- and then only very large
coefficients are detected.

A specific case of interest is for
$$\tau_n=O(\sqrt{\log{p}/n}).$$
For instance, it can easily be shown that this case especially
occurs with overwhelming probability when the entries $X$ are i.i.d. gaussian $N(0,\sigma^2)$. This case is considered in details in
the following simulation study.
\begin{corollary}\label{consistence_gaus}
Assume $\tau_n=O(\sqrt{\log{p}/n}).$
 The algorithm
LOL$(X,Y,O(\sqrt{\log{p}/n}),O(\sqrt{\log{n}/n}))$ is consistent as
soon as the vector $\beta$ satisfies to the sparsity, size and
significance conditions for
$$\mu_n\geq O(\sqrt{\log{p}\log{n}/n})\quad\mbox{ and }\quad
S\leq O(\log{n}).$$
\end{corollary}

\section{Adaptive LOL and  practical properties}\label{section adaption}
In this part we address in an experimental way some problems which
are not solved in the theoretical part. Among them, the most crucial
question is the way to choose the tuning parameters $t_n$ and $s_n$.
We focus on building an auto driven procedure which is illustrated with
some simulations.

\subsection{LOLA Procedure}
The auto-driven procedure  is called
  LOLA procedure, which is LOL( ) completed by procedure A( ) allowing the adaptive tuning of the
  threshold levels $t_n$ and $s_n$ \label{LOLTH}

$$
\mbox{LOLA}(X,Y)=\mbox{LOL}(X,Y,\hat t,\hat s)
$$
with the following  choices of the tuning parameters
\begin{eqnarray*}
\hat t= \mbox{A}((\tilde{X_{1}}^tY,\ldots,\tilde{X_{p}}^tY))&\mbox{
and } &\hat s=\mbox{A}((\widehat{\beta_1},\ldots,\widehat{\beta_p}))
\end{eqnarray*}
where $\tilde{X}$ and $\hat \beta$ are defined in LOL() algorithm.
The algorithm A( ) is described by

\vspace{0.5cm}

\begin{center}
\label{algoth}
\begin{tabular}{|l|ll|}\hline
{\bf $u\leftarrow$ A$(Z)$}&
{\bf Input:} & variables $Z=(Z_1,\ldots,Z_m)$\\
&{\bf Output:} &level $u$\\ \hline\end{tabular}\end{center}

\vspace{0.5cm}

\noindent and the output $u$ is computed as follows. Let
$|Z|_{(1)}\leq |Z|_{(2)}\leq \ldots \leq |Z|_{(m)}$ be the ordered
 sample and consider the deviance function defined by
 $$
 \mbox{dev}(J)=\sum_{j=1}^J\left(|Z|_{(j)}-
 \overline{|Z|}^{(J-)}\right)^2+\sum_{j=J+1}^m\left(|Z|_{(j)}-
 \overline{|Z|}^{(J+)}\right) ^2
 $$
 where $ \overline{|Z|}^{(J-)}$ and $ \overline{|Z|}^{(J+)}$ are
 the empirical means of the $|Z|_{(j)}$'s for
 respectively $j=1,\ldots,J$ and $j=J+1,\ldots,m$. We choose as
 threshold level
 $$u = |Z|_{(\widehat{ J})}\quad \mbox{ for }\quad
\hat J=\mbox{Arg}\min_{j=1,\ldots m}\mbox{dev}(J).
$$

Notice that A() is entirely data-driven and can be roughly justified
as follows: since the
 thresholds are used to select the higher responses $|Z|'$s (being here
 either the scalar product $|\tilde{X}^tY|'$s or the estimators of the
 linear coefficient $|\hat \beta|'$s) among a set of given variables,
 we aim at splitting the set of
variables,
 into two clusters in such a way that the higher ones are forming one of the two
 clusters. The output of the A( ) algorithm is the frontier computed between the
 clusters: it is then chosen by minimizing the deviance
 between classes (see also Kerkyacharian, Mougeot, Picard and Tribouley (2009) ).
Obviously, A( ) algorithm performs better when both clusters are
distinctly separated which is the case in our theoretical setting since the
sparsity assumption suggests that the law of the $|Z|$ (in absolute
value) should be a mixture of two distributions: one for the
variables included in the model (positive mean) and one
 for the others which should be very small (zero mean).

Observe that the A( ) algorithm has its own interest since it could
be useful in many other
 nonparametric settings such as denoising or density
estimation where a thresholding procedure is performed. For
instance, the input of A( ) could be the empirical wavelet
coefficients $\hat \beta_{j,k}$ when local thresholding is
considered.

\subsection{Illustration with simulations}\label{illu}The performances of
 LOLA algorithm are first presented by considering a classical
 framework
where the predictors are realizations of gaussian variables.
Intensive studies have been performed and we present here the
results for $n=400$ and $p=2000$. LOLA procedure is repeated $K=100$
times using each time $n$ different  random observations.
 Observations  are simulated from the model $Y=X\beta+\epsilon$ where $\epsilon$ is a
gaussian vector such that the  signal over noise  ratio satisfies
$SNR=5$. The cardinal of the set of indices of the predictors to be
selected $\cI=\{\ell=1,\ldots,p,\,\beta_\ell\not=0\}$ is denoted
$S$. Three experiments are presented
\begin{itemize}
\item Exp1: the predictors $X_1,\ldots,X_p$ are independent and
$S=10$
\item Exp2: the predictors $X_1,\ldots,X_p$ are linear dependent and
$S=10$
\item Exp3: the predictors $X_1,\ldots,X_p$ are independent
and $S=50$.
\end{itemize}

The empirical coherences computed on the $n-$sample of the
predictors $X=(X_1,\ldots,X_p)$ are

\vspace{0.2cm}

\begin{center}
\begin{tabular}{rrrr}
  \hline
 &Exp1&Exp2&Exp3\\
 $\tau_n$&0.25&0.72&0.25\\
  \hline
\end{tabular}
\end{center}

\vspace{0.2cm}

\noindent and deliver the following message: the results for Exp2
should be considered with cautiousness since $\tau_n$ is large. It
is also a benefit of LOLA to compute an empirical indicator giving a
warning sign to the user.

\begin{figure}[h!]
 \begin{center}
 \begin{tabular}{ccc}
\includegraphics[height=4cm]{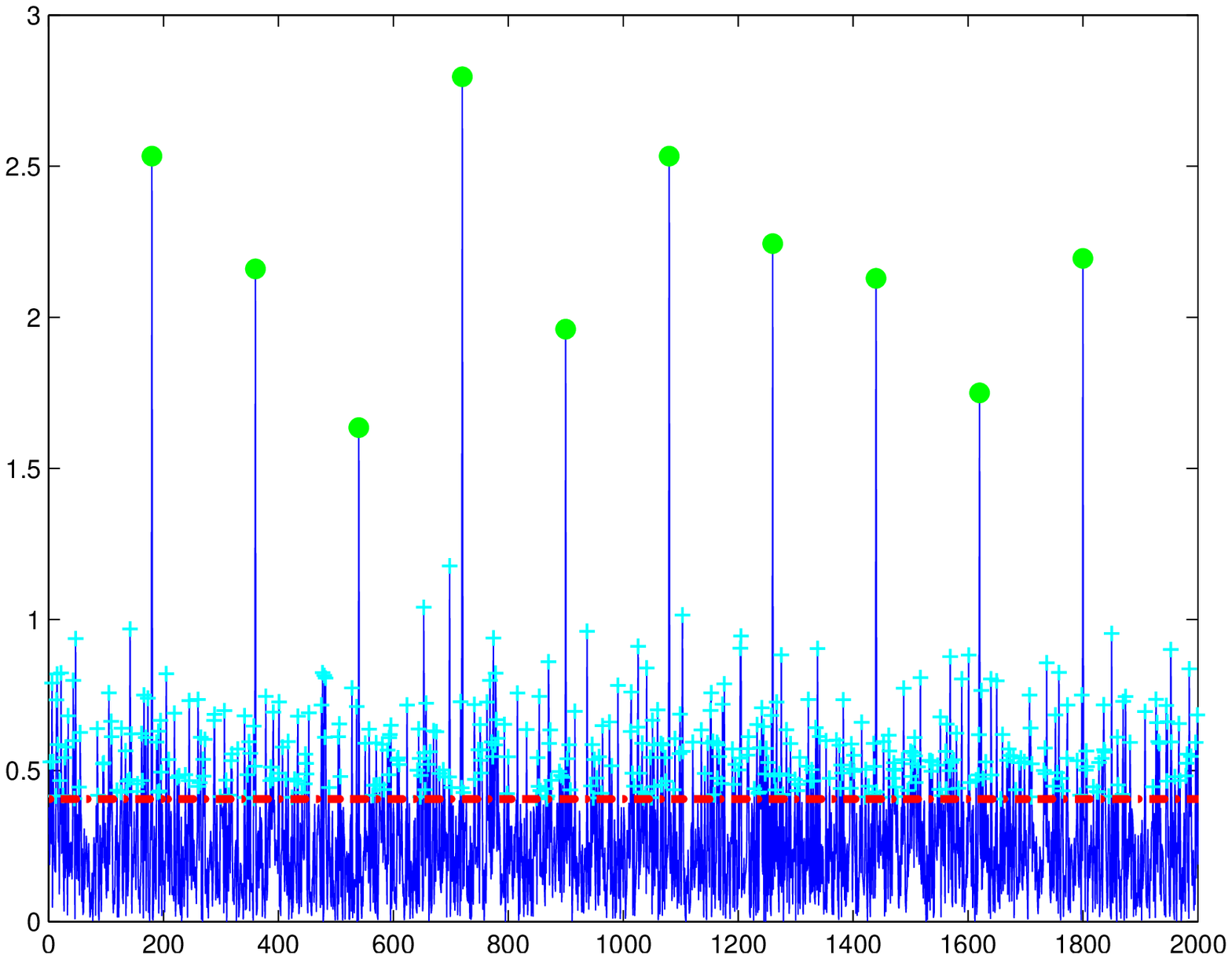} &
\includegraphics[height=4cm]{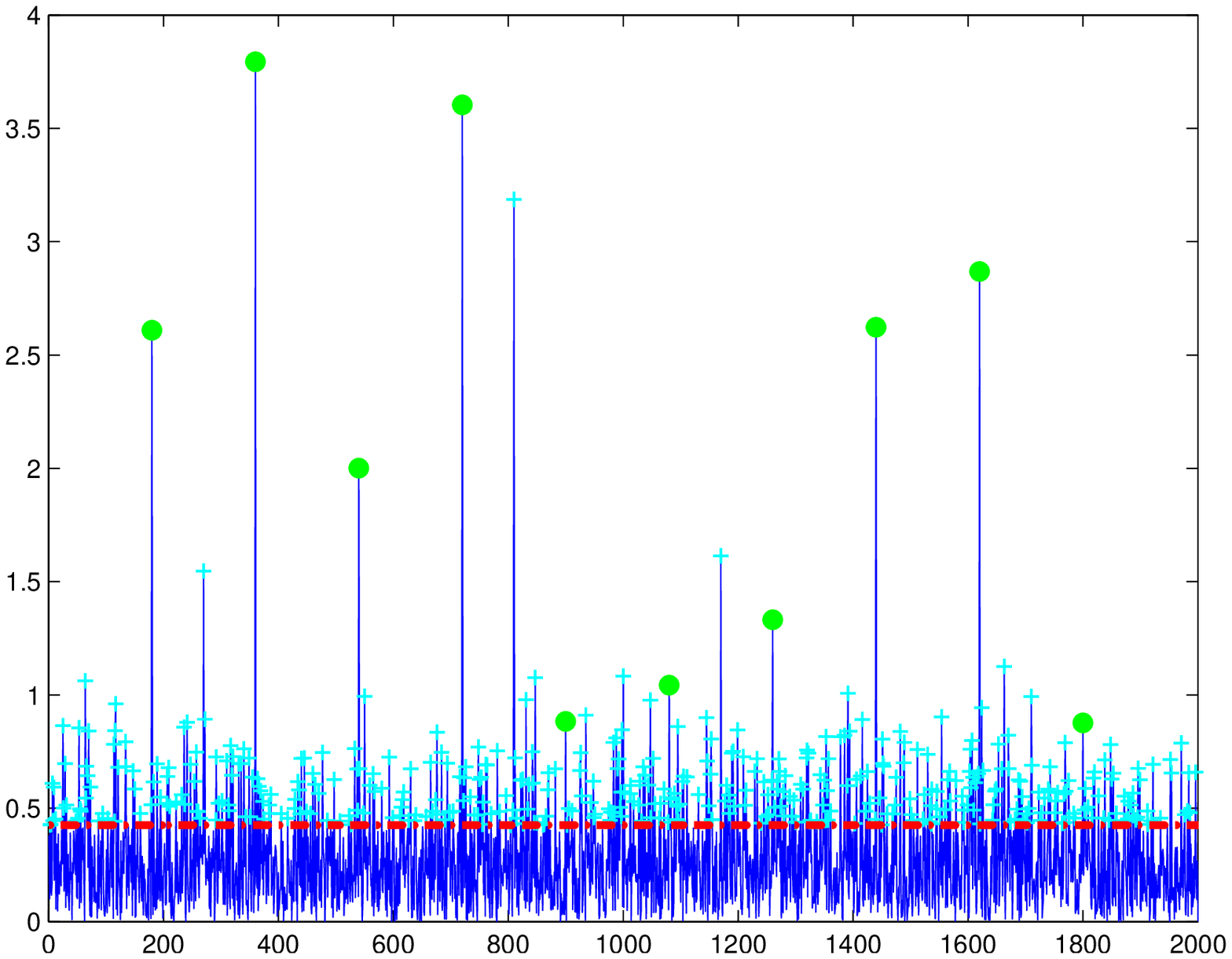} &
\includegraphics[height=4cm]{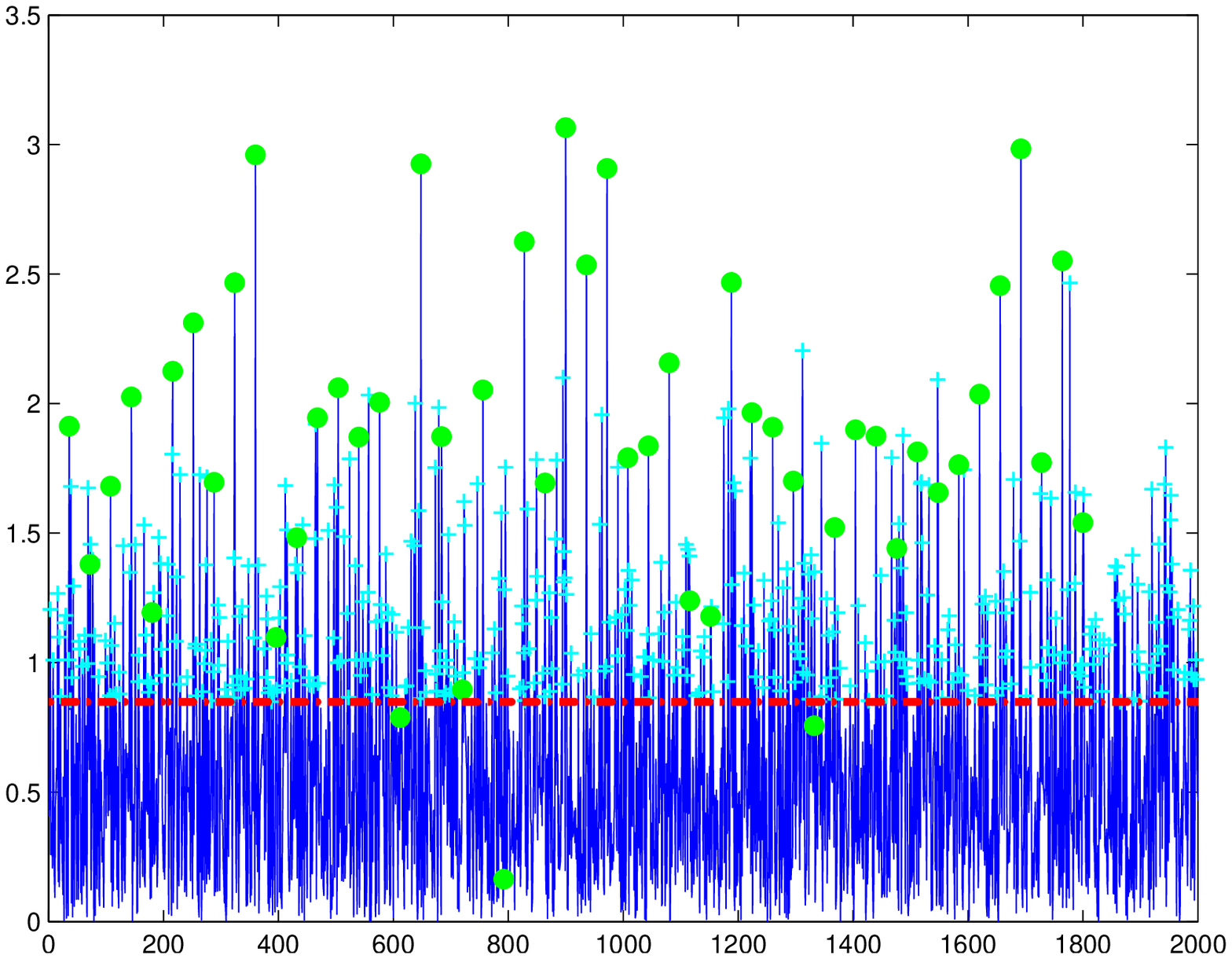}\\
Exp1&Exp2&Exp3
 \end{tabular}
 \end{center}
\caption{Empirical second moment between $Y$ and $X_1,\ldots,X_p$
for $p=2000$, $n=400$, the predictors $X'$s are gaussian. Exp1: The
predictors are independent, $S=10$. Exp2: The predictors are
dependent, $s=10$. Exp3: The predictors are independent, $S=50$. The
horizontal line is the auto-driven level $\hat t$.}
 \label{fig:1}
\end{figure}

The two successive thresholding steps of the algorithm are first
detailed using the three different examples presented above. Figure
\ref{fig:1} illustrates the first step of LOLA for the selection of
the leaders for one experiment  chosen among all $K=100$
experiments. All the scalar products  $|\widetilde{X_\ell}^tY|$ for
$\ell=1,\ldots p=2000$ are represented in the picture; the estimated
level $\hat t$ computed with  procedure $A( )$ is indicated with a
horizontal line. The leaders are all the variables
$\widetilde{X_\ell}$ with a scalar product
 $|\widetilde{X_\ell}^tY|$ exceeding the threshold
 and are labeled with a small cross. let $\cL$ denotes the
 set of leader indices. The indices
belonging to $\cI$ are circled to indicate the variables which should be
rightly selected.

 Observe that the reduction of the dimensionality is
drastic: $N=144$ leaders are selected during the first step among
the $p=2000$ initial predictors. When the sparsity is small and when
the predictors are independent (see Exp1), the values of the scalar
products of the predictors really involved in the model are close to
the value of the coefficients (here $|\beta_\ell|\sim 2$). For Exp1,
the first step is fine since any variable defined in $\cI$ is chosen
to be a leader. As the sparsity increases or when the predictors become
dependent as it is the case for Exp2 or Exp3, the empirical scalar products between the
predictors and the target may have more unstable behaviors leading to significant coefficients falling below the
 threshold $\hat t$ and  as a consequence not selected as leaders during
 the first thresholding step. For Exp3,
three variables $p_0\in\{612,790,1338\}$ which should be kept
are eliminated at the first step and are definitively lost for selection.

\begin{figure}[h!]
 \begin{center}
 \begin{tabular}{ccc}
\includegraphics[height=4cm]{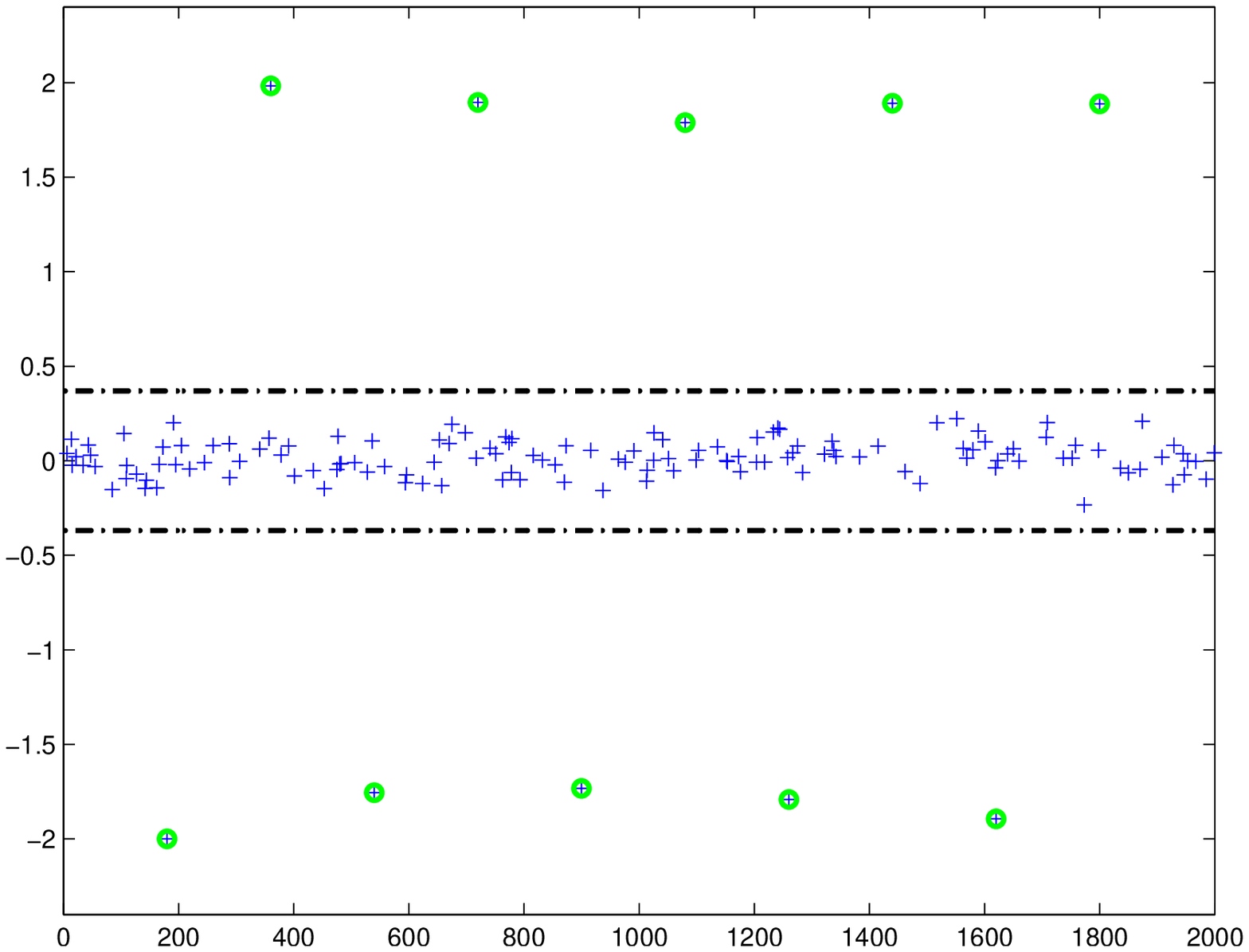} &
\includegraphics[height=4cm]{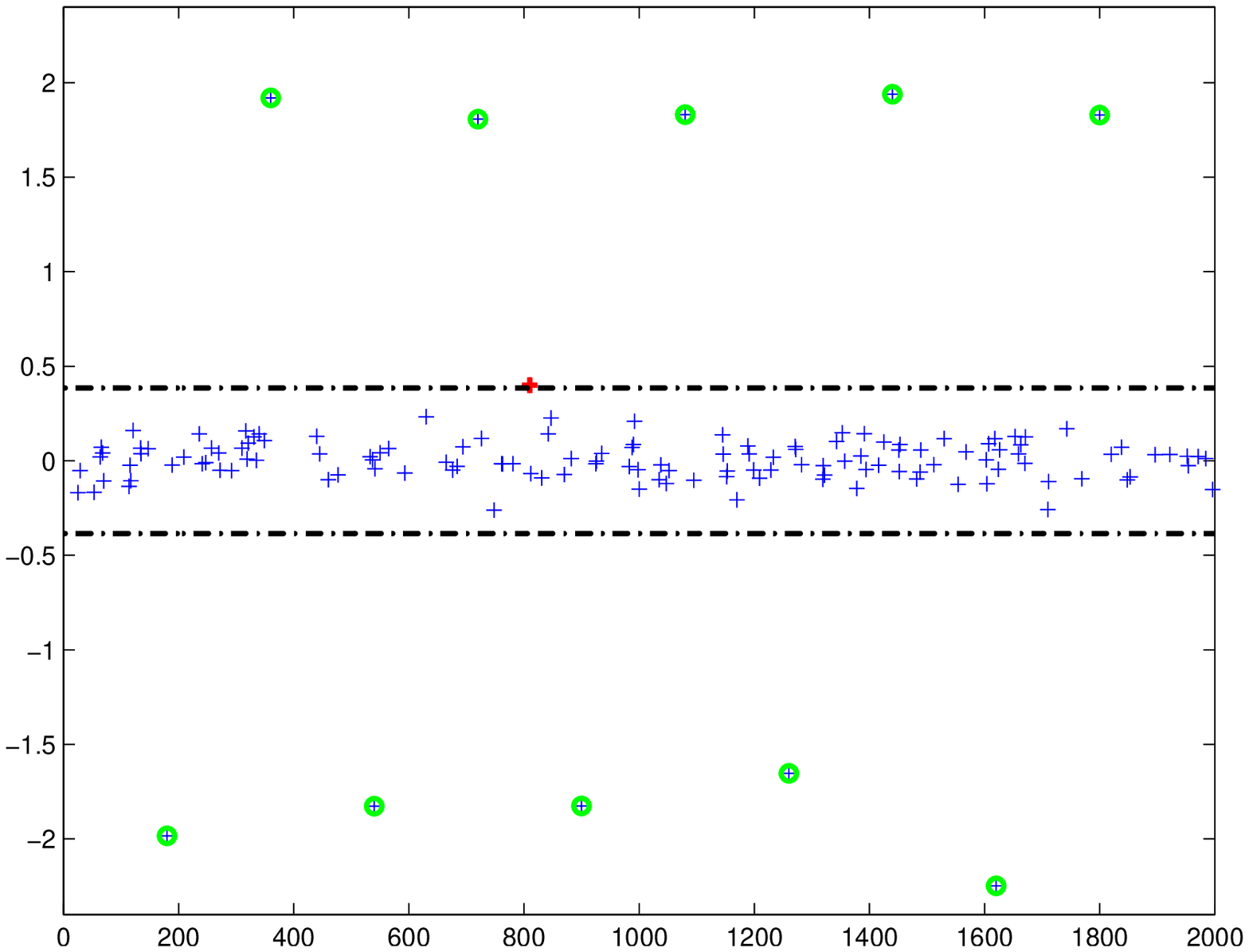} &
\includegraphics[height=4cm]{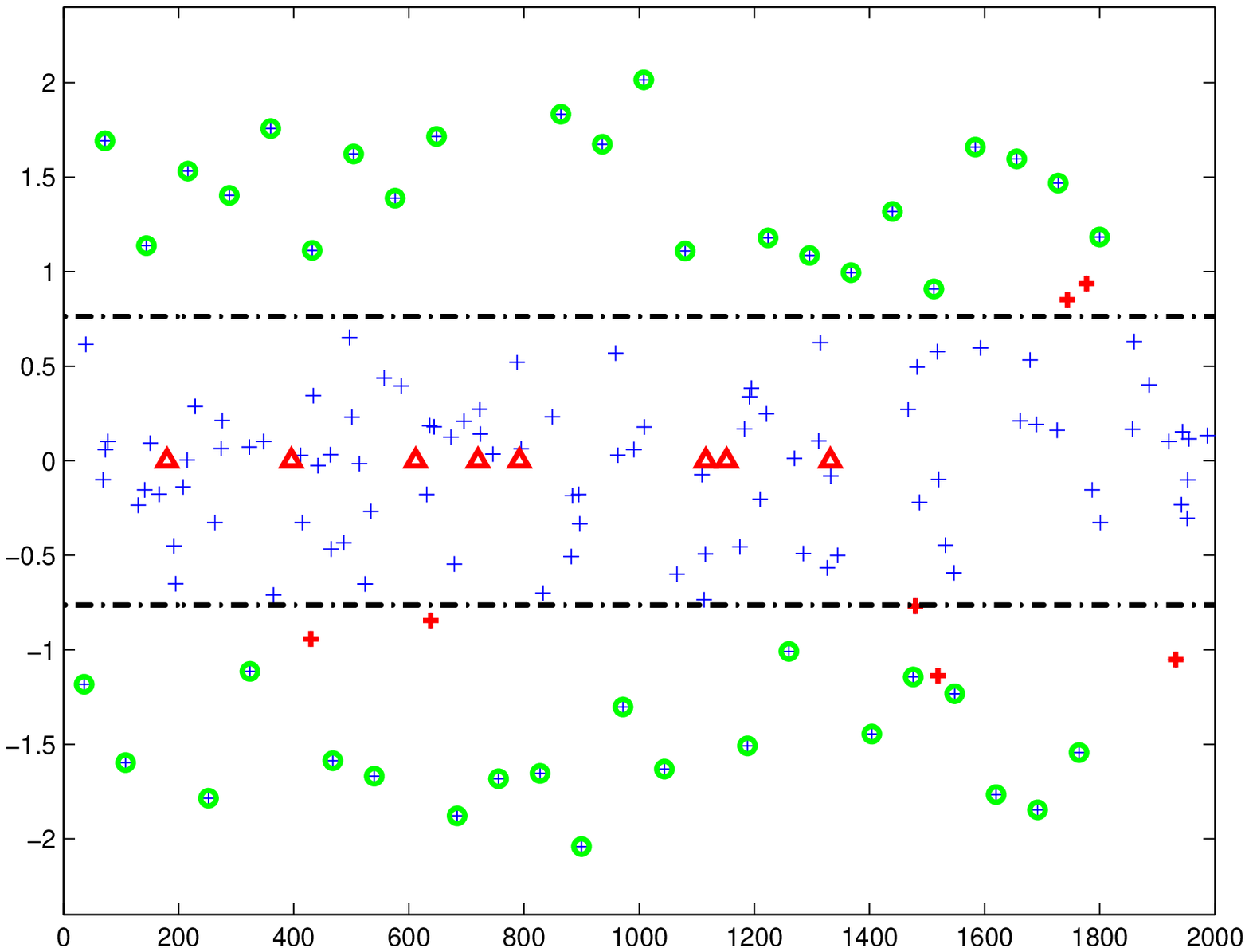}\\
Exp1&Exp2&Exp3
 \end{tabular}
 \end{center}
\caption{Estimators of the linear coefficients $\beta_\ell$ for
$\ell\in\cL$ for $p=2000$, $n=400$, the predictors $X'$s are
gaussian. Exp1: The predictors are independent, $S=10$. Exp2: The
predictors are dependent, $S=10$. Exp3: The predictors are
independent, $S=50$. The horizontal lines are the auto-driven levels
$\pm\,\hat s$.}
 \label{fig:2}
\end{figure}

Figure \ref{fig:2} illustrates the effect of the second thresholding
of LOLA procedure for the same experiments for Figure \ref{fig:1}.

The horizontal lines represent the levels
 $\pm\hat s$ where $\hat s$ is the output of procedure A( ). Figure \ref{fig:2}
 gives
 the estimated coefficients
$\widehat{\beta_\ell}$ for $\ell\in\cL$ computed with the OLS method
 restricted to the sub-set of predictors selected as
leaders. The circle is the label for the coefficients
$\widehat{\beta_\ell}$ with $\ell\in \cI$ which are kept: this label
allows to see the number TP of true positive detections. Triangular
is the label for the coefficients $\widehat{\beta_\ell}$ with
$\ell\in \cI$ but which are below the level $\hat s$ (in absolute
value): this label allows to see the number FN of false negative
detections. Diamond is the label for the coefficients
$\widehat{\beta_\ell}$ with $\ell\not\in \cI$ which are kept by
LOLA: this label allows to see the number FP of false positive
detections. The results for Exp1 are excellent: both clusters of
coefficients
 $\widehat{\beta_\ell}$ are well separated and  A( ) algorithm
 performs in this situation  very accurately; we obtain FN=FP=$0$ and we get exactly
 $\widehat{\cI}=\cI$.  For Exp3, the separation between both
 clusters is not so straight and miss detections (triangular pattern) as false detections
 (cross not circled) are observed:
 FN$=8$ and FP$=7$.

The global performances obtained with $K$ experiments  are presented
in Table 1. For each column, the number of indices to be estimated
is indicated into brackets. The first two columns focus on true
detections and the last ones on false detections. Results for exp1
are perfect: LOLA algorithm selects always the right predictors and
there is no error. When the sparsity $S$ increases, the results are
less powerful: some predictors (FN$=13.9$) which should be selected
by LOLA are finally not selected. Nevertheless, we observe that the
number of false positive detections is again small (FP$=1.2$).

\subsection{Conclusion and Comments}

In order to be concise in this presentation, the above illustrations
focused only on gaussian random variables  with $SNR=5$. It should
be stressed that an extensive simulation study has been conducted in
parallel in order to evaluate the performances of LOLA. Various
distributions for the predictors and different SNR have been
implemented and studied. As conclusions to this experimental design,
it can be stressed that for a given number of observations and
potential predictors, the selection is more accurate for a low
sparsity level. The number of false negatives and false positives
tend to be null as the number of observations increases for a fixed
number of predictors and/or the sparsity level decreases for a fixed
number of predictors. It should be noted that very similar results are
obtained if the Gaussian predictors are replaced by uniform, Bernoulli random variables or mixture of the above distributions.

 As can be seen in Corollary \ref{consistence_gaus},
LOL procedure is consistent under a condition on the coherence
$\tau_n\le O\left(\sqrt{\log{ p}/n}\right)$.
 This condition is
verified with overwhelming probability for instance when the entries
of the matrix $X$ are independent and identically random
variables with a sub-gaussian common distribution but the results
obtained in Exp2 show that the procedure is still working quite well even if
this hypothesis is not satisfied. This fits a common fact in high dimensional setting, which is that the theoretical results are often more pessimistic than the true performances of the procedures.
Moreover, before running the
procedure, it should be noticed that the computation of $\tau_n$
brings some benefit as an indication of potential misbehaviors.

\section{LOLA properties in a toy 'instrumental' setting}\label{Boston}

We begin by studying the practical quality of our algorithm with
real data combined with  simulated instrumental data by revisiting
the Boston Housing data set available from the {\it UCI machine
learning data base repository}: http://archive.ics.ucfi.edu/ml/. One
of our goal is here to evaluate  the ability of LOLA procedure  to
accurately select the right predictors even if the original
variables are embedded in a huge space built with artificial
variables. This analysis also help to point out
 which kind of variables are wrongly selected from the complementary artificial
 space and which one are not selected from the initial space.

 The original Boston Housing data are defined by one continuous
target variable $Y$ (the median value of owner-occupied homes in
USD) and $p_0=13$ predictive variables observed for $n=506$
observations. For our purpose, these original data are embedded in a
high dimensional space of size $p=2100+13$ by adding artificially
$300$ independent random variables of $7$ different laws: normal,
$\log$normal, bernoulli, uniform, exponential, Student and Cauchy in
equal proportion. This set of distributions is chosen to mimic the
different empirical laws of the $p_0=13$ original predictive
variables. In order to numerically evaluate the performances of
LOLA, the procedure is called $K=100$ different times. Each time, a
new set of artificial 'instruments' has been simulated. Also, to
evaluate the instability of the algorithm, each time the procedure
has been performed on a 0.75$n$ sub-sample of the initial $n$-
sample, randomly chosen.

Since the initial data are observed (and not simulated), we do not
know in advance which variables should be selected in the model. In
order to evaluate the performances of  LOLA algorithm, we considered two benchmark
procedures: a) the classical multiple regression following by
variable selection using a simple Student test, b) the stepwise
regression method (significant level 95\%). Obviously, the two benchmark procedures
are performed   in the regular space of data with $p_0=13$ variables and $K=100$ times again by randomly choosing  and $0.75\,n$ observations among the original data set.

In the high dimensional model ($p=2113$), the empirical coherence is
$\tau_n=0.98$ which is very high and  indicates that the
predictors are very dependent. We applied LOLA procedure to this huge
set of data. The results show that the number of false detections of
the artificial variables is extremely low: we only select $92$
adding variables over  2100 random variables for a total of $K=100$
experiments. It is interesting to observe that half of the
 selected variables are distributed according to
the cauchy law and $10\%$ are distributed according to the $T(2)$
law. A complementary work shows that the impact of heavy tailed
distributions for the predictors
 is similar to the impact of dependence between predictors.

Figure \ref{fig:3} shows the frequencies of detection for the
initial $p_0=13$ predictors using LOLA procedure, OLS with Student
test and stepwise regression. The results obtained with LOLA and OLS
with Student test are similar very similar. This comparison confirms
that our procedure performs fairly well in the presence of a huge
number of artificial 'instrumental variables'.

This preliminary investigation, as a first step, as well as the
previous simulation study, justifies the use in the sequel of LOLA
as a selection algorithm in the presence of an important set of
instrumental variables and make stronger our conclusions of the
following section.


\section{ International Economic Growth}\label{PNB}

 We study in this section the problem of convergence hypothesis in
 economic expounded in the introduction. We use the Barro and Lee data available from
http://www.nber.org/pub/barro.lee. We present different models to
evaluate the casual effect of a initial level $X$ of gross domestic
product on the growth rate (of gross domestic product) using
parametric models and economic variables as well as non-parametric
models. Our aim is to empirically prove that the convergence
hypothesis is valid and to compare our results with the results
presented in Belloni and Chernozhukov (2010).

The Barro and Lee data are extracted as in Belloni and Chernozhukov
(2010) and missing data are removed for each studied case. Notice
that
 we restrict the number of countries ($n$) in the study since we want to take into consideration a large number of
variables ($p$). It is important to keep this in mind for a further interpretation
of the results.

\subsection{The data}
The national growth rate in gross domestic product (GDP) per capita
is our dependent variable $Y$ and is studied for different periods of time.
The predictor $X$ is an initial amount of the gross domestic product
(GDP) per capita. More precisely, the variables are defined by
$Y=\log\left({\mbox{GPD}_{t_2}/\mbox{GPD}_{t_1}}\right)$ and
$X=\log(\mbox{GPD}_{t_0})$ for the following periods of time:

\vspace{0.5cm}

\begin{center}\begin{tabular}{ccccc}
\hline &Period $[t_1,t_2]$&Initial Date $t_0$& $n$ & $p$
\\\hline
Exp1&1965-75&1960&63&208\\
Exp2&1965-75&1965&63&208\\
Exp3&1975-85&1970&52&375\\
Exp4&1975-85&1975&52&375\\
\hline
\end{tabular}
\end{center}

\vspace{0.5cm}

The Barro and Lee data contain different economical indicators from
1960 to 1985. Six categories of variables are considered: Education
(1), Population/Fertility (2), Governement Expendidures (3), PPP
deflators (4), Political variables (5), Trade Policy and others (6).

In view to measure the accuracy of LOLA, as well as to increase the
stability of the algorithm and also to see the impact of the sample of
countries on the selection and estimation, we run the procedure
$K=1000$ times using each time on a portion of $0.85n$ data randomly
chosen from the initial data set. All the given results (the
estimator $\hat \alpha_2$, the bounds of the confidence interval
computed under the gaussian hypothesis, the $R^2$ and the p-value
associated to the global Fisher test) are then averaged for the
$K=1000$ experiments. The empirical standard deviation normalized by
$\sqrt{K}$ is given into brackets. The indicator $N_0$ is the
frequency of the even "zero belongs to the confidence interval".
Ideally, $N_0=0$.
 The confidence intervals are computed for a coverage of $90\%$
 which leads to test the null hypothesis $H_0:\;\alpha_2=0$ against $H_1:\;\alpha_2<0$ at the level $5\%$.

\subsection{LOLA in a parametric setting}\label{PNBlola}
First, it should be noticed thanks to the results obtained by the
standard OLS given in Table 2 that
 the linear model $Y=\alpha_1+\alpha_2X+u$ is irrelevant. For all
 periods of time,  the $R^2$ is almost zero and the   p-value leads to reject the significance of the model.
   Since $N_0=1$, the hypothesis $H_0:\;\alpha_2=0$ is always accepted at level $5\%$.

We now use LOLA procedure  to select the subset
 of explanatory variables $Z_{\mbox{selected}}$
 containing the maximum of information in order
to explain $X$ or $Y$. Taking inspiration into the vast literature
about instrumental variables in econometrics (see among many others
Angrist, Imbens and Rubin (1996), Blundell and Powell (2003),
Florens (2003), Darolles, Fan, Florens, Renault (2010) and Florens,
Heckman, Meghir, Vytlacil (2003) in a nonparametric or
semi-parametric framework), we select the instruments using the two
following models
$$
\mbox{Model 1:}\quad X=Z\beta+\mbox{error}\quad\mbox{ and }\quad
\mbox{Model 2:}\quad Y=\alpha_1+\alpha_2X+Z\beta+\mbox{error}.
$$
Depending on the situation, each of them has its own justification
and interest: Model 1 is used in the first step of the Instrumental
Variable method while Model 2 relies the endogeneity of $X$ in the
initial model $Y=\alpha_1+\alpha_2X+u$ with the possibility that
covariates $Z$ are missing.
 The selected instruments $
Z_{\mbox{selected}}$ (obtained either via Model 1 or Model 2) are
then used as control variables and we consider the following  model
\begin{eqnarray}\label{modelcontrole}
Y=\alpha_1+\alpha_2X+ Z_{\mbox{selected}}\beta+\mbox{error},
\end{eqnarray}
where $\alpha_1,\; \alpha_2$ and $\beta$ are estimated (OLS).
 For Model 1 (above) and Model 2 (below), we obtain the following results given in Table 3.

The selection via Model 2 seems to be more appropriate than the
selection via Model 1: the p-value associated to the Fisher test is
very small which leads to accept the significance of Model (Equation
\ref{modelcontrole}). We now focus on the results obtained with the
selection using Model 2. Observe that the results concerning Exp1
and Exp2 (respectively Exp3 and Exp4) are very similar: the date
$t_0$ of the initial amount of the GDP leads to similar results. The
estimation of the coefficient $\alpha_2$ is negative and the
confidence interval never contains zero (remember that the results
are averaged on $K=1000$ experiments). Results are better for Exp1
(and Exp2) than for Exp3 (and Exp4): $N_0$ is equals $2\%$ for the
first period of study instead of $30\%$ for the second period of
study. Analyzing the empirical densities of the estimator $\hat
\alpha_2$ given in Figure 4, we notice that for both periods of
time, the support is almost included in $R^-$: for the period
$1965-1975$ and $t_0=1960$, all of the runs provide negative
estimated values while for the period $1975-1985$ and $t_0=1970$,
99.8\% of the runs do.

Notice that the number  of selected instruments is quite reasonable ($\hat
S\sim 6$ or $8$) so we can also comment our results on a qualitative
point of view. Incidentally, it is very interesting  to identify the
determinants of growth. To partially answer  this question, Figure
5 shows the frequencies of selected variables using
Model 2 when the growth rates under consideration are for both
periods 1965-75 and 1975-85. The six  vertical areas define the
6 broad categories of variables as given in the beginning of this
part.
It is interesting to notice that the sets of selected
variables are quite different for both periods of time. For 1965-75, the most selected
variables are indicators of the demography of the countries like
"Life expectancy at age 0" (1960, 1970, 1965-69, 1970-74) and
"Total Fertility rate" (1965-69)
while for period 1975-85, they are
"Total gross enrollment ratio for secondary education"  (1975),
"Male gross enrollment ratio for secondary education" (1965),
"Growth rate of population" (1965-69),
"Ratio of real government "consumption" expenditure net of spending on defense and on education to real GDP" (1965-69) and
"Black market premium" (1975-79, 1980-84).
This, of course can be partially due to the instability of the method,
it may also have economical interpretations.

\subsection{Comparison with results by Belloni and Chernozhukov}
In Belloni and Chernozhukov (2010), a lasso procedure is used for
the first
 period 1965-75. The method depends on a tuning
 parameter $\lambda$: the smaller is $\lambda$,
 the more instrument variables are selected and larger
  (in absolute value) is the estimator of $\alpha_2$.
  In the following table, we recall the results by
  Belloni and Chernozhukov (2010) when $\lambda$ is varying and
  our results when all the data are used ($K=1$)
  since no stability bootstrapping has been performed in
  Belloni and Chernozhukov (2010).
  Observe that the LOLA procedure leads to larger estimators of
    $\alpha_2$ (in absolute value) and confidence intervals which
    are better separated from zero. We verify that zero never belongs to the intervals
 of confidence obtained by LOLA procedure even if the level $\alpha$
 becomes larger than $90\%$.

For 1965-75, the selected variables are "Percentage of "no schooling
in the male population" (1970) and "Percentage of "primary school
attained" in the total population and in the female population"
(1970). These variables related to the education policy are also
selected by Belloni and Chernozhukov (2010). It is striking that
LOLA selects variables which are the same as in Belloni and
Chernozhukov (2010). Recall that LOLA works in high dimension (here
$p\sim 200$ or $300$) while the lasso procedure is used in the
standard case where $p<n$ ($p\sim 40$). For the second period
1975-1985, LOLA retains
\begin{itemize}
\item "Percentage of primary school complete in the female population"  (1960, 1965)
\item "Population Proportion over 65" (1960)
\item "Growth rate of population" (average on 1980-84)
\item "Black market premium" (average on 1975-79, average on 1980-84)
\end{itemize}
as in Belloni and Chernozhukov (2010) and in addition
\begin{itemize}
\item "Total gross enrollment ratio for secondary education" (1975, 1980)
\item "Total fertility rate" (1970, average on 1980-84)
\item "Life expectancy at age 0" (1975).
\end{itemize}

\subsection{LOLA in a non-parametric setting}
If we are only interested by the sign of the coefficient $\alpha_2$,
it could be interesting to build a set of tailored instrumental
variables which have no economical meaning but which provides a very
good model
$$
Y=\alpha_2X+\alpha_1+Z\beta+\mbox{ error}
$$
in the sense that the number $\hat S$ of variables $Z$ is small and
that the null hypothesis $H_0: \alpha_2=0$ is rejected. Using our own instruments has also another advantage: the study can be conducted with almost all the countries ($n\sim 110$) while  half were removed in the previous part because pieces of information were missing for some variables: by instance China, Singapore were not considered for the first period of time and many countries from Africa were missing in the study for the second period.

As a consequence, let us consider the vector $X$ as a signal, and
indeed  replace $X$ by $X_{(.)}$ where $X_{(.)}$ is the vector such
that $X_{(1)}\leq \ldots\leq X_{(n)}$ such a way that the curve of
the initial amount of the GDP becomes smoother. Taking inspiration
from  the learning theory (see Kerkyacharian, Mougeot, Picard and
Tribouley (2009) ), we build  a set of functions called dictionary
${\mathcal D}$ containing the functions cosinus, sinus, box and
Schauder:
$$
{\mathcal D}=\left\{\phi_\lambda\right\}_{\lambda\in\Lambda}=
\left\{\phi^1_{\lambda_1},\phi^2_{\lambda_2},\phi^3_{\lambda_3,}\phi^4_{j,k}\quad\mbox{
for } 1\leq\lambda_1,\lambda_2,\lambda_3\leq n,
0\leq j\leq\log_2(n),1\leq k\leq 2^j-1 \right\}
$$
for
\begin{eqnarray*}
\phi^1_\lambda(x)&=&\cos(2\pi\lambda\,x)\quad\mbox{ and }\quad
\phi^2_\lambda(x)=\sin(2\pi\lambda\,x)\\
\phi^3_\lambda(x)&=&1_{[a_\lambda,b_\lambda]}(x)
 \mbox{ where } a_\lambda,b_\lambda \backsim
{\mathcal U}_{[0,n]}\\
\phi^4_{j,k}(x)&=&2^{j/2}\phi^3(2^jx-k)\quad\mbox{ where }\quad\phi^4(x)=\left(x1_{[0,0.5]}(x)\right)
-\left(x1_{[0.5,1]}(x)-1\right).
\end{eqnarray*}
For each function of the dictionary, the curve $Z_\lambda=\;^t\left(\phi_\lambda\left(1/n\right),
\ldots,\phi_\lambda\left(i/n\right),\ldots,\phi_\lambda\left(n/n\right)\right)$  is considered as an instrumental variable and the matrix $Z=\left(Z_\lambda\right)_{\lambda\in\Lambda}$ contains $p=\#{\mathcal D}$ predictors. Here, we get $n=110$ and $p=300$. The obtained results are excellent

Observe that the number of selected instruments is smaller than
in the parametric setting as well as the  fit is
better if we are only concerned by the explanation
of $Y$ by $X$. The confidence interval for $\alpha_2$ almost
never contain zero: see the very small values of $N_0$. Again, we observe
that the initial $t_0$ has no influence on the estimation of $\alpha_2$ and that
the estimator is slightly smaller for the second period of time. In order to compare the non
parametrical methodology with the parametrical methodology, we give again the
 empirical densities of the estimator $\hat \alpha_2$: here the results are
 excellent since their supports are strictly contained in $R^-$.

Finally, if we end this section with a rapid comparison with Belloni
and Chernozhukov (2010), using all the available data in one run, we
always accept the hypothesis $\alpha_2<0$. Notice that the number of
selected instruments is again relatively small and is comparable to
the number founded in Belloni, Chernozhukov (2010)
 for $\lambda=0.4$. Again, our estimators are further away from zero than the
 lasso estimators.

\section{Proof of Theorem \ref{consistence_MS}}\label{section_proof}

Theorem \ref{consistence_MS} is an obvious consequence of
Proposition \ref{false_detection}.

%
The main ingredients for  proving Proposition \ref{false_detection},
are Lemma 1 and Lemma 2. These lemmas are proved in Mougeot, Picard
and Tribouley (2010): see the terms $IBS$ and $OBB+OBS$ for Lemma 1
and the term $IBS$ for Lemma 2 and put
$\alpha_\ell=\sigma_\ell\,\beta_\ell$ and
$\widehat{\alpha_\ell}=\sigma_\ell\,\widehat{\beta_\ell}$.
\begin{lemma}Recall that $\cJ$ is the set of the indices of the leaders defined in (\ref{leader}). Assume that the sparsity conditions
 (\ref{Scond}) are satisfied by the sequence
 $(\beta_\ell),\ell=1,\ldots,p$.
 For any
$\lambda\geq C\left(S\tau_n^2+S\log{p}/n\right)$ where $C$ is a
computable constant depending on $\eta^2$ and $\sum_{\ell=1}^p|\sigma_\ell\beta_\ell|$, we
get
\begin{align*}
P\left(\sum_{\ell\in\Omega_\iota}\;|\sigma_\ell\beta_\ell|^2\;\;>\;\lambda\right)
\leq 2\exp\left(-c\,\frac{n\lambda}{\eta^2}\right)\quad\mbox{ for
}\iota=1,2,
\end{align*}
with
$$\Omega_1=\{\ell=1,\ldots,\;  p, \;|\sigma_\ell\beta_\ell| \geq
2s_n,\,|\sigma_\ell\widehat{\beta_\ell}| \leq s_n\}$$
and
$$
\Omega_2=\{\ell=1,\ldots,\; p,\;\ell\not\in\cJ,\,|\sigma_\ell\beta_\ell|
\geq 2t_n\}$$ and where $c$ is an universal constant.
\end{lemma}

\begin{lemma} Recall that $\cJ$ is the set of the indices of the leaders defined in
(\ref{leader}). Assume that the sparsity conditions
 (\ref{Scond}) are satisfied by the sequence
 $(\beta_\ell),\ell=1,\ldots,p$ and choose threshold levels $t_n,s_n$
 such that
$t_n>s_n\vee C\left(\tau_n+\sqrt{\log{p}/n}\right)$. Then, for any
$\lambda\geq C\left(S/n\right)$ where $C$ is a constant depending on
$\sigma^2$ and $M$, we get
\begin{align*}
P\left(\sum_{\ell\in
\Omega_3}\;|\sigma_\ell\widehat{\beta_\ell}-\sigma_\ell\beta_\ell|^2\;
\;>\;\lambda\right)\leq
\exp\left(-c\,\frac{n\lambda}{\sigma^2}\right)
\end{align*}
with
$$\Omega_3=\{\ell\in\cJ,\;|\sigma_\ell\beta_\ell| \leq
s_n/2,\,|\sigma_\ell\widehat{\beta_\ell}| \geq s_n\}$$ and where $c$ is an
universal constant.
\end{lemma}
Now, we are ready to prove Proposition \ref{false_detection}.
First, we focus on $FN$.
Using \eref{B2cond}, and for $\mu_n>s_n/2$, we get
  $$ \#(\cI)=\#\left(\{\ell,\;|\sigma_\ell\beta_\ell|\geq \mu_n\}\right)\le S\le \lfloor\tau_n^{-1}\rfloor$$
and we conclude applying Lemma 1 with $\mu_n\geq 2(s_n\vee t_n)$.
More precisely, for a given $k$, we get
\begin{align*}
P(FN>k)&=P(\#(\cI\cap (\widehat{\cI})^c) \;>\;k)\\
&=P(\#\{\ell,\;|\sigma_\ell\beta_\ell|\not =0,|\sigma_\ell\widehat{\beta_\ell}|\,
\I\{|\sigma_\ell\widehat{\beta_\ell}|\geq s_n\}=0\}\;>\;k)\\
&\leq
P(\sum_{\ell=1}^{p}|\sigma_\ell\beta_\ell|^2\;\I\{\,|\sigma_\ell\beta_\ell|
\geq \mu_n\}\;\I\{\,|\sigma_\ell\widehat{\beta_\ell}|\,
\I\{|\sigma_\ell\widehat{\beta_\ell}|\geq s_n\}=0\}
\;>\;k\mu_n^2)\\
&\leq
P(\sum_{\ell\in\Omega_1}\;|\sigma_\ell\beta_\ell|^2\;>\;k\mu_n^2/2)+
P(\sum_{\ell\in\Omega_2}\;|\sigma_\ell\beta_\ell|^2\;>\;k\mu_n^2/2)\\
& \leq 2\exp\left(-c\;\frac{nk\mu_n^2}{2\sigma^2}\right)
\end{align*}
as soon as
$$
k\;\mu_n^2/(2\sigma^2) \geq c^\prime\,\left(M^2\;\tau_n^2 \vee \frac
1n\vee \frac{|\log{\tau_n}|}{n}\right)\,S.
$$
Similarly, applying Lemma 2 with $t_n>s_n\vee
C\left(\tau_n+\sqrt{\log{p}/n}\right)$, we get
\begin{align*}
P(FP>k)
&=P(\#(\cI^c\cap \widehat{\cI}) \;>\;k)\\
&=P(\#\{\ell,\;|\sigma_\ell\beta_\ell| \ge \mu_n,|\sigma_\ell\widehat{\beta_\ell}|\,
\I\{|\sigma_\ell\widehat{\beta_\ell}|\geq s_n\}\not =0\}\;>\;k)\\
&\leq
P\left(\sum_{\ell\in\Omega_3}\;|\sigma_\ell\widehat{\beta_\ell}-\sigma_\ell\beta_\ell|^2\;
\;>\;ks_n^2\right)\;\leq
\exp\left(-c\;\frac{nks_n^2}{\sigma^2}\right)
\end{align*}
a soon as
$$
k\;\frac{s_n^2}{\sigma^2} \geq c\,\frac Sn.
$$

\vspace{0.5cm}

 \noindent {\Large \bf References}

\vspace{0.3cm}

{\small
\begin{description}
\item[ {\rm Angrist, J., G. Imbens, and D. R. Rubin.}] (1996).
Identification of Causal Effects Using Instrumental Variables.
\textit{Journal of the American Statistical Association},
\textbf{91}, pp 444--455.

\item[{\rm Barro, R. and  Lee, J.-W.}] (1994). {\sl Data set for a
panel of 139 countries.}

\item[ {\rm Barron A. R., Cohen A., Dahmen, W. and DeVore, R. A.}] (2008).
Approximation and learning by greedy algorithms. \textit{ The Annals
of Statistics}, \textbf{36}, pp 64--94.

\item[{\rm Belloni, A. and Chernozhukov, V}] (1999).  $L_1$-penalized quantile regression in
high-dimensonal sparse models.  \textit{ ArXiv e-prints}

\item[{\rm Bickel, P. J., Ritov, Y. and Tsbybakov, A.}] (2008).
Simultaneous analysis of Lasso and Dantzig selector. \textit{ ArXiv
e-prints}

\item[{\rm Blundell, R. and Powel, J.}] (2003).
Endogeneity in Nonparametric and Semiparametric Regression Models,
 \textit{Cambridge Univ. Press.},
{Advances in Economics and Econometrics: Theory and Applications},
\textbf{2}, pp 312--350.

\item[{\rm Bunea, F., Tsbybakov, A. and Wegkamp, M.}] (2007).
Sparsity oracle inequalities for the Lasso,
 \textit{Electronic Journal of Statistics}, \textbf{1}, pp 169--194.

\item[{\rm Candes, E. and Tao, T..}] (2007).
The Dantzig selector: statistical estimation when $p$ is much larger than $n$
 \textit{The Annals of Statistics}, \textbf{35}, pp 2313--2351.

\item[{\rm Darolles, S. and Fan, Y. and Florens, J.-P. and Renault, R.}] (2010).
 Nonparametric Instrumental Regression, \textit{IDEI Working Paper}, \textbf{228
}.

\item[{\rm Fan, J. and Lv, J.}] (2008).
A Selective Overview of Variable Selection in High Dimensional
Feature Space.
 \textit{ Journal of the Royal Statistic Society B}, \textbf{70},
pp 849--911.

\item[{\rm Fan, J. and Lv, J.}] (2010).
Sure independence screening for ultrahigh dimensional feature space
 \textit{ Statistica Sinica}, \textbf{20},
pp 101--148.

\item[{\rm Florens, J.P..}] (2003). Inverse problems and structural
 econometrics: The example of instrumental variables.
 \textit{Cambridge Univ. Press., In Advances in Economics and Econometrics: Theory and Applications},
\textbf{2}, pp 284-311.

\item[{\rm Florens, J.-P. and Heckman, J. and Meghir, C. and Vytlacil, E.}] (2003).
Instrumental Variables, Local Instrumental Variables and Control
Functions,  \textit{IDEI Working Paper }, \textbf{249}.

\item[{\rm Haury, A.-C., Jacob, L. and Vert, J.-P.}] (2010).
Increasing stability and interpretability of gene expression
signatures, \textit{ArXiv e-prints }.

\item[{\rm
Kerkyacharian, G., Mougeot, M., Picard, D. \& Tribouley, K.} ]
(2009). Learning Out of Leaders, \textit{Multiscale, Nonlinear and
Adaptive Approximation.} Lecture Notes in Comput. Sci. Springer.

\item[ {\rm  Mougeot, M., Picard, D. \& Tribouley, K.}] (2010). Learning Out of Leaders.
\textit{ArXiv e-prints}

\item[ {\rm Needell, D. and Tropp, J. A.}] (2009).
Co{S}a{MP}: iterative signal recovery from incomplete and inaccurate
samples \textit{Applied and Computational Harmonic Analysis},
\textbf{26}, pp 301--321.

\item[{\rm Sala-I-Martin, X.}] (1997). I Just Ran Two Million
Regressions. \textit{The American Economic Review}, \textbf{87}, pp
178--183.

\item[ {\rm Tibshirani, R.}] (1996). Regression Shrinkage and Selection via the
Lasso. \textit{Journal of the Royal Statistical Society B},
\textit{}, \textbf{58}, pp 267-288.

\item[ {\rm Tropp, Joel A. and Gilbert, Anna C. }] (2007).
Signal recovery from random measurements via orthogonal matching
pursuit. \textit{Institute of Electrical and Electronics Engineers.
Transactions on Information Theory}, \textbf{53}, pp 4655--4666.

\end{description}}

\newpage

\newpage

\newpage

\begin{center}
\begin{table}[h!]
\begin{tabular}{rrrrr}
  \hline
 &  TN & TP & FN & FP \\
  \hline
Exp1& 1990.0 (1990) & 10.0 (10) & 0.0 (0) & 0.0 (0)\\
Exp2& 1989.9 (1990) & 9.2 (10)  & 0.8 (0) & 0.7 (0)\\
Exp3& 1948.8 (1950) & 36.1 (50) & 13.9 (0)& 1.2 (0) \\
  \hline
\end{tabular}
\caption{ Simulations performances for $K=100$ replications. The
true target value is given into brackets.\label{tableFNFP}}
  \end{table}

\end{center}

\newpage

\begin{figure} [h!]
 \begin{center}
 \begin{tabular}{ccc}
\includegraphics[width=5cm,height=3cm]{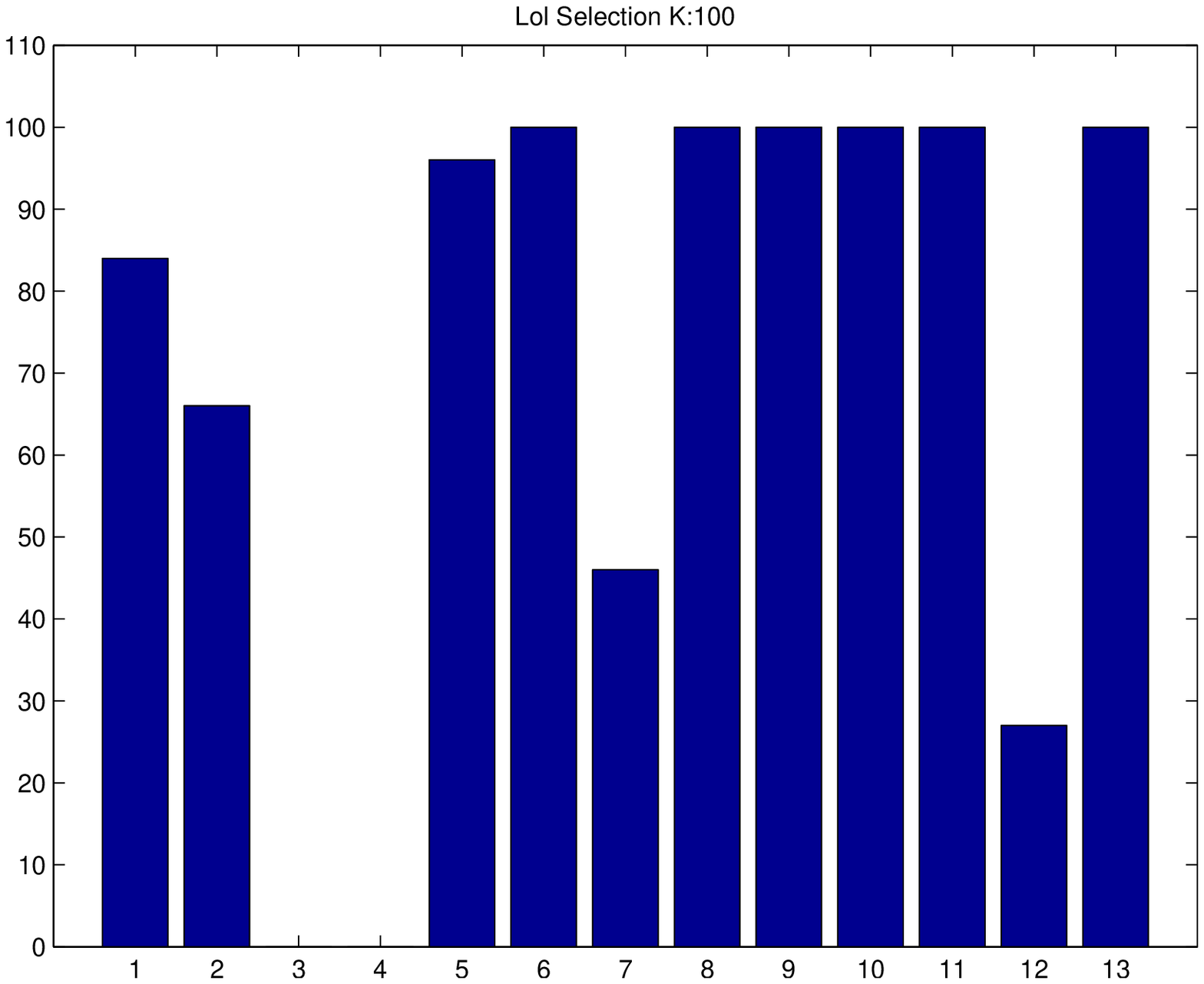}&
\includegraphics[width=5cm,height=3cm]{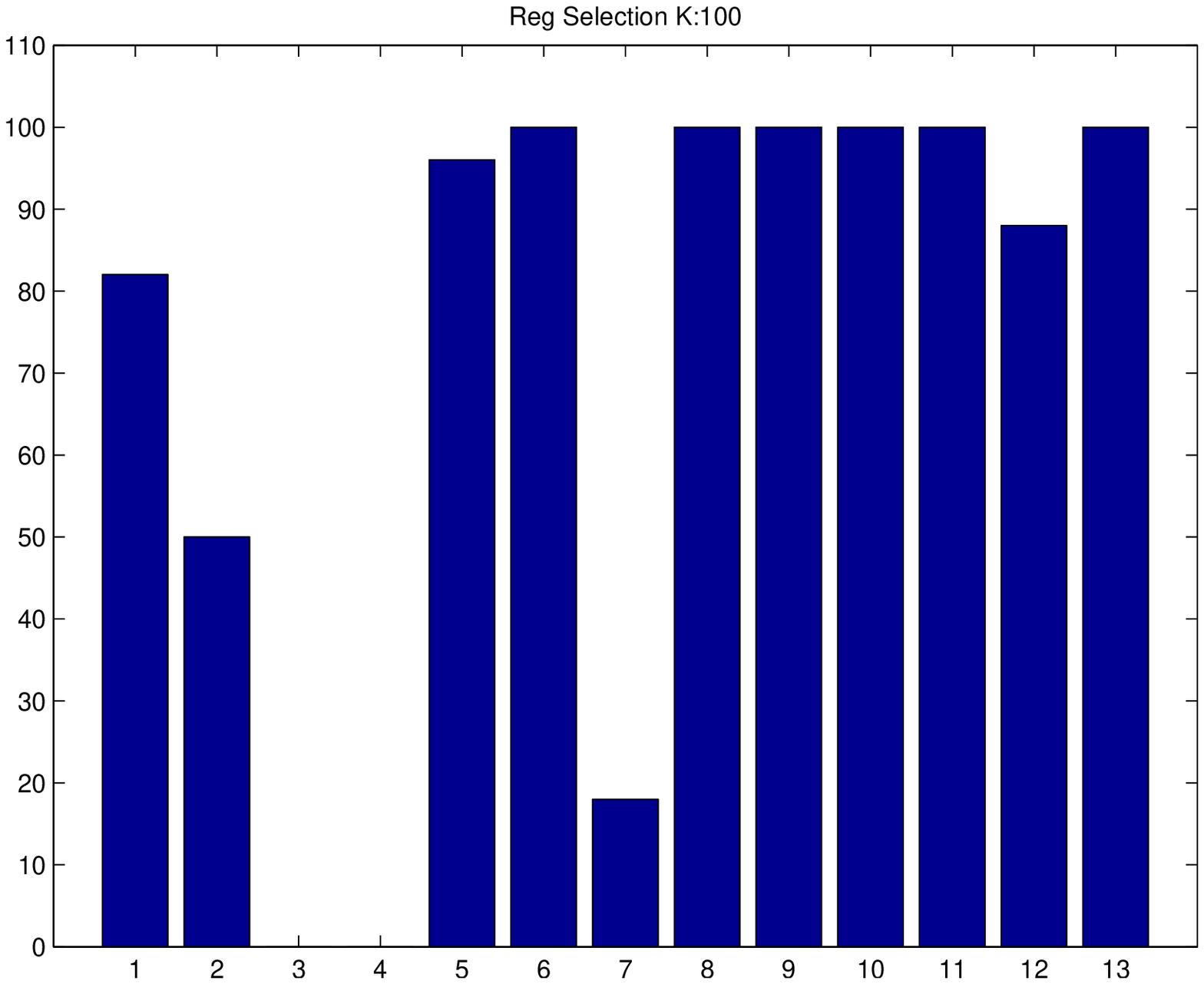}&
\includegraphics[width=5cm,height=3cm] {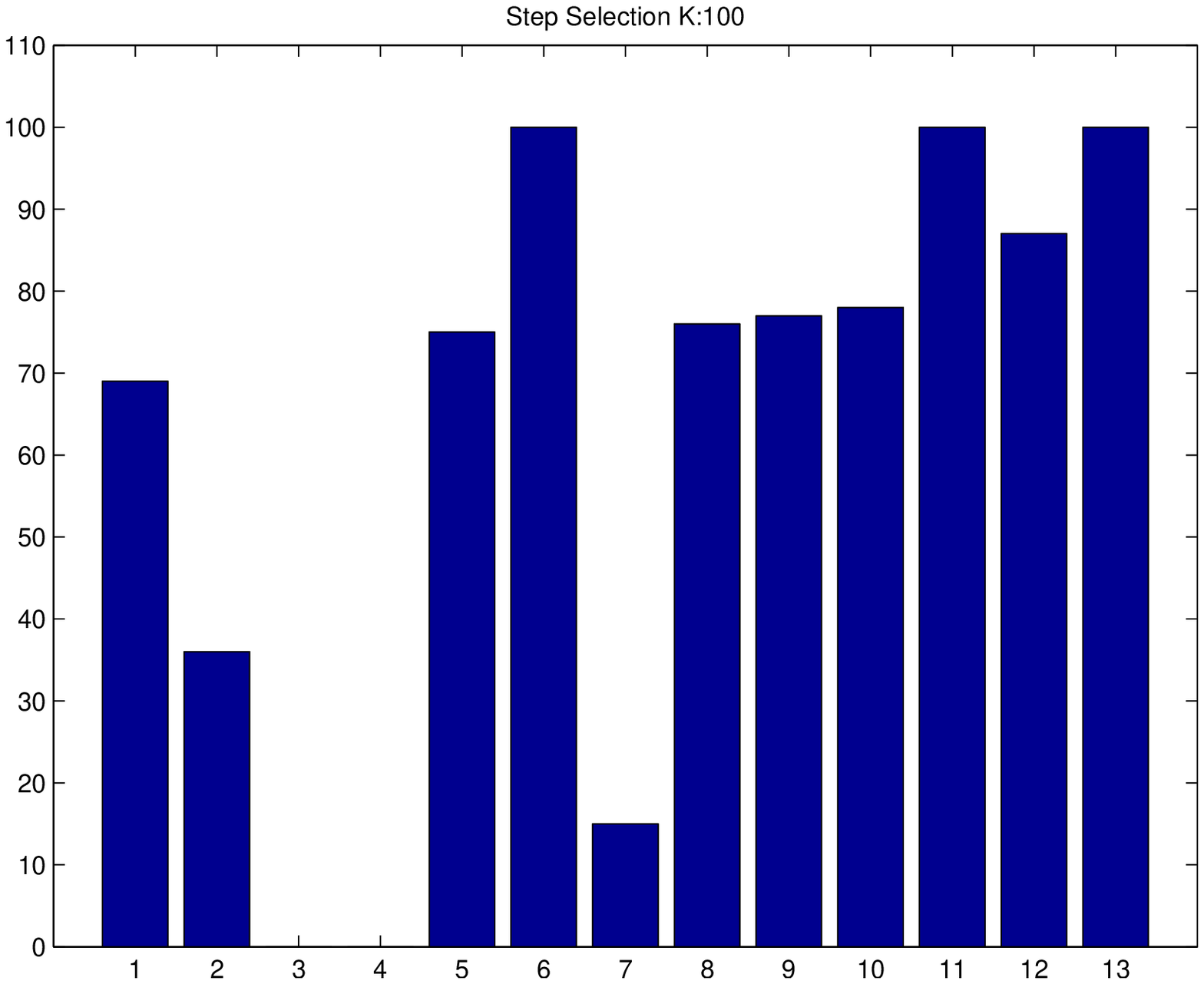}\\
LOLA ($p=1213$) & OLS Student Selection($p=13$)& STEP ($p=13$)\\
\end{tabular}
\caption{Variables selection performances for Boston data, $K=100$
\label{BostonPerf}}
 \end{center}
 \end{figure}

\newpage

 \vspace{0.5cm}
 {\footnotesize
\begin{center}
\begin{tabular}{cccccc}
\hline &$\widehat{\alpha_2}$&Confidence Interval &$N_0$& $R^2$&
p-value for $F$ test
\\\hline
Exp1&0.012  & [-0.038 , 0.062  ] & 1.000 & 0.005  & 0.673  \\
Exp2&0.017  & [-0.031 , 0.065 ] & 1.000 & 0.008  & 0.577   \\
Exp3&0.051  & [-0.002 , 0.104 ] & 1.000 & 0.045  & 0.139   \\
Exp4&0.058  & [\,0.006 , 0.111  ] & 1.000 & 0.060  & 0.085  \\\hline
\end{tabular}
 \end{center}
{TABLE 2: Estimation of the model $Y=\alpha_1+\alpha_2X+u$ for the
different periods of time. All the standard deviation values
normalized by $\sqrt{1000}$ are smaller than $10^{-3}$.}}

\newpage

\vspace{0.5cm}

\begin{center}{\footnotesize
\begin{tabular}{rccccc}
\hline &$\hat S$&$\widehat{\alpha_2}$&CI &$N_0$& p-value for $F$
test
\\\hline
Exp1&10.8 (0.1) & -0.245 (0.002) &[-0.358 (0.002), -0.132 (0.002)] & 0.011 (0.003) & 0.013 (0.001)\\
Exp2&10.9 (0.1) & -0.225 (0.002) &[-0.342 (0.002), -0.109 (0.002)] & 0.039 (0.006) & 0.025 (0.002)\\
Exp3&14.2 (0.1) & -0.227 (0.004) &[-0.429 (0.004), -0.025 (0.004)] & 0.452 (0.015) & 0.155 (0.005)\\
Exp4&14.6 (0.1) & -0.175 (0.004) &[-0.434 (0.005),  0.083 (0.004)] & 0.742 (0.014) & 0.255  (0.006)\\
&&&&\\
Exp1&6.3 (0.1) & -0.246 (0.002) &[-0.342 (0.002), -0.150 (0.001)] & 0.016 (0.004) & 0.007 (0.041)\\ 
Exp2&6.3 (0.1) & -0.232 (0.002) &[-0.331 (0.002), -0.132 (0.001)] & 0.025 (0.005) & 0.011 (0.042)\\ 
Exp3&8.4 (0.1) & -0.166 (0.003) &[-0.266 (0.003), -0.066 (0.002)] & 0.230 (0.013) & 0.004 (0.017)\\ 
Exp4&8.4 (0.1) & -0.149 (0.003) &[-0.262 (0.003), -0.036 (0.002)] & 0.311 (0.014) & 0.009 (0.034)\\ 
\hline
\end{tabular}

{TABLE 3: Estimation of $\alpha_2$ for determinants $Z$ selected via
Model 1 (above) and Model 2 (below). }}
\end{center}

\newpage

 \begin{center}
 \begin{tabular}{cc}
\includegraphics[height=5cm,width=8cm]{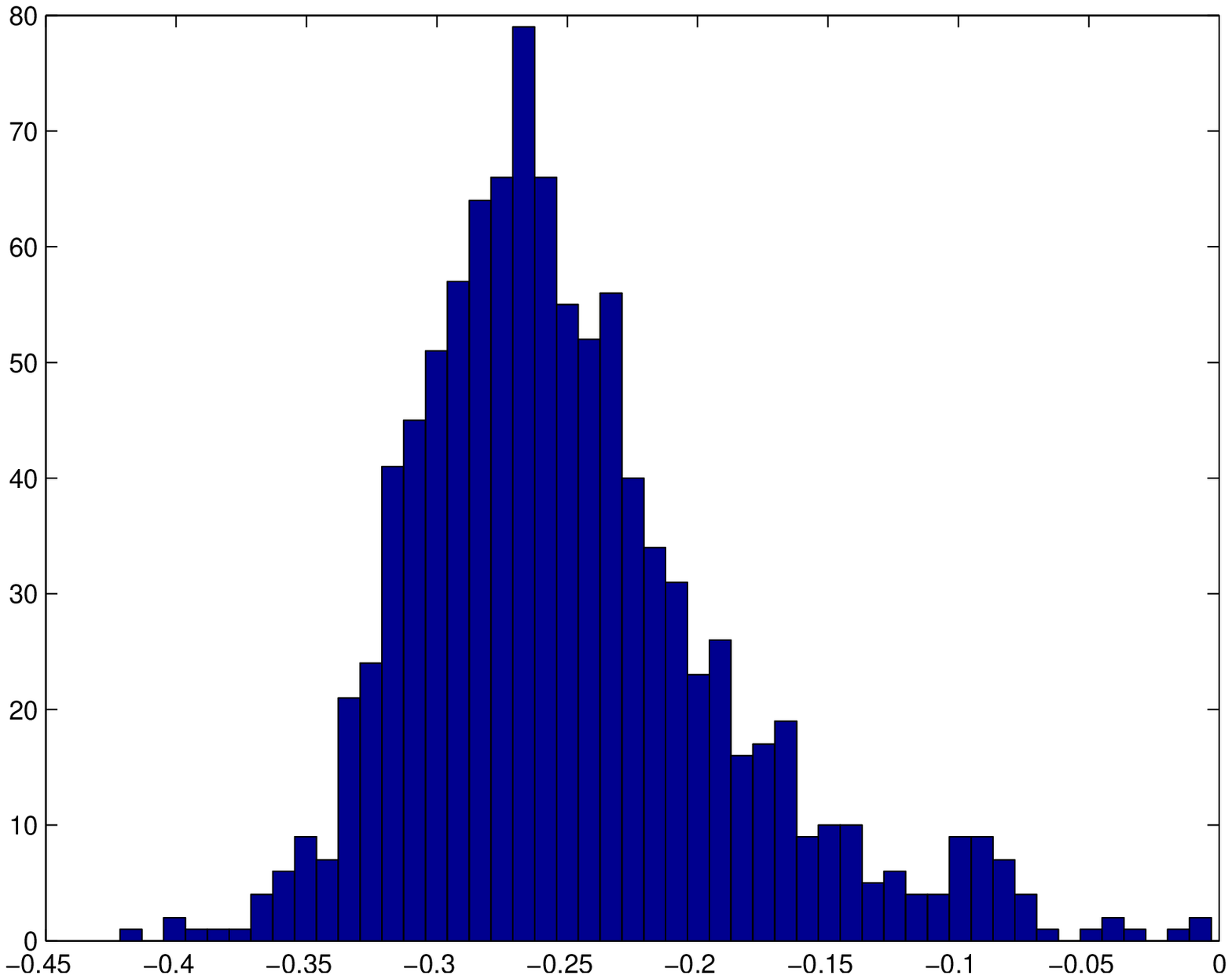} &
 \includegraphics[height=5cm,width=8cm]{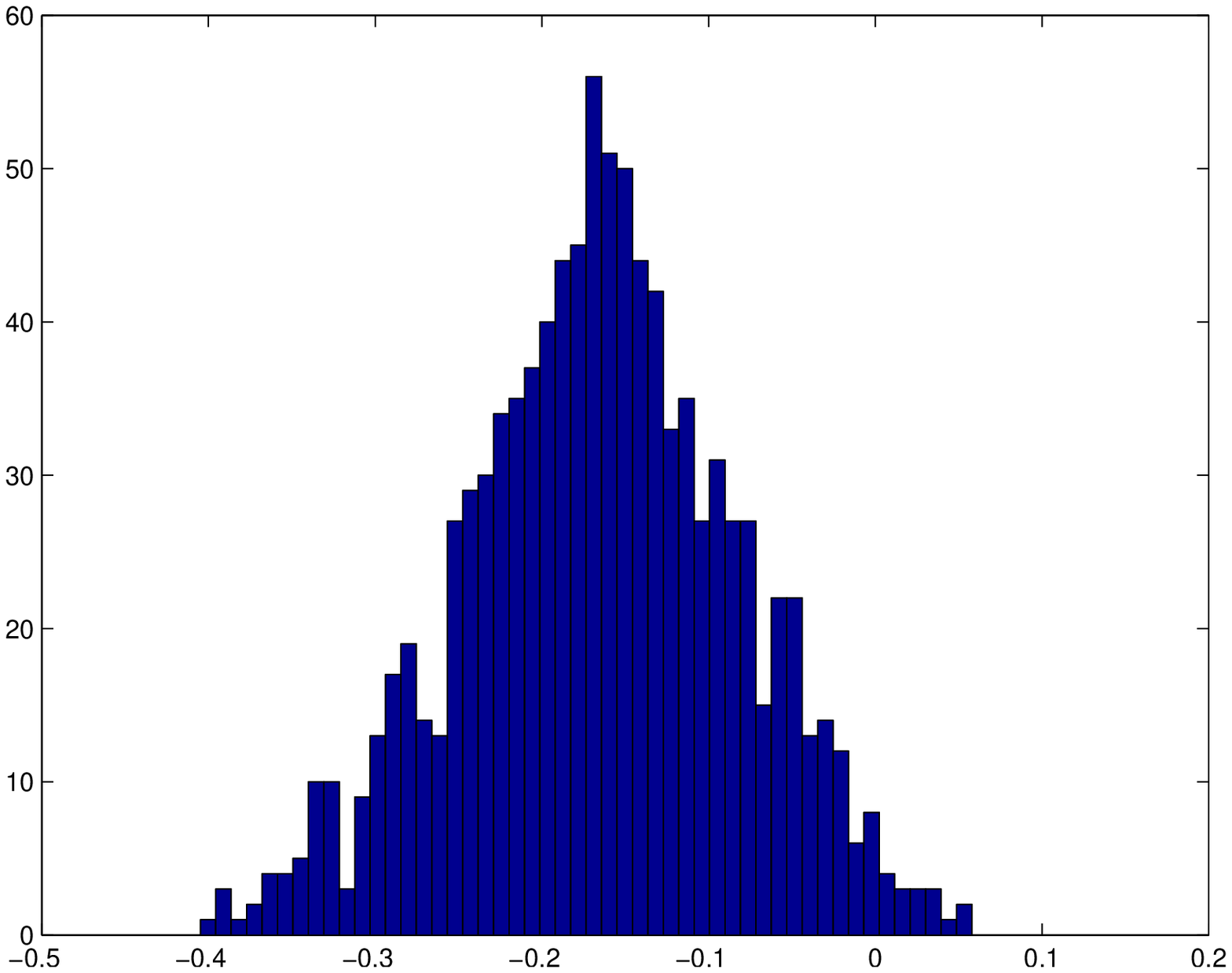}\\
{\footnotesize FIGURE 4. Exp1: 1965-75, $t_0=1960$, empirical
density of $\hat\alpha_2$}&{\footnotesize Exp3: 1975-85, $t_0=1970$,
empirical density of $\hat\alpha_2$}
 \end{tabular}
 \end{center}

\newpage

 \begin{center}
 \begin{tabular}{c}
\includegraphics[height=5cm,width=12cm]{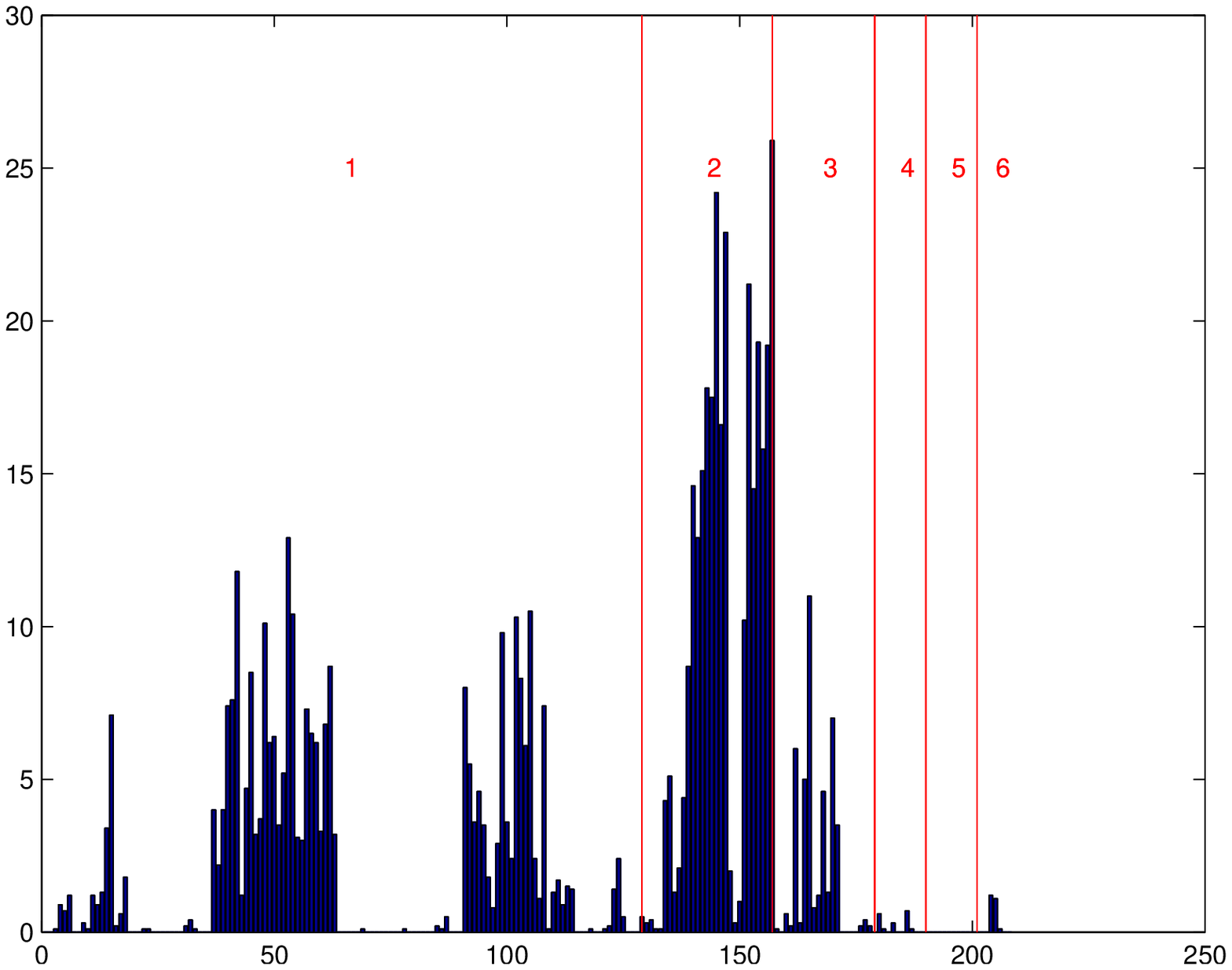}
\\
\includegraphics[height=5cm,width=12cm]{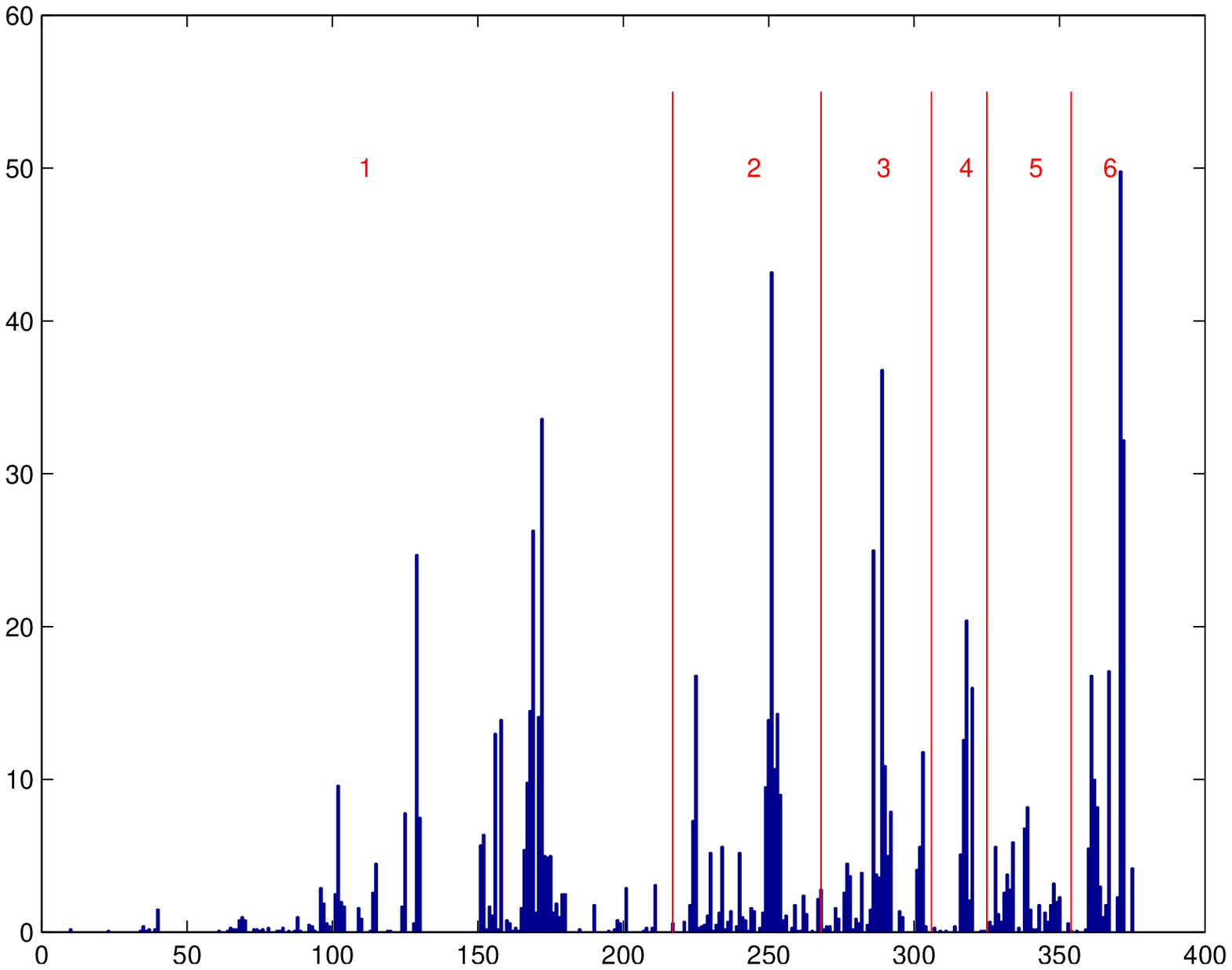}
 \end{tabular}
 \end{center}
{\footnotesize{FIGURE 5. Selected variables for different periods
1965-1975, $t_0=$ 1960 (bottom) and 1975-1985, $t_0=$1970 (below).
Area 1: Education, Area 2: Population/Fertility, Area 3: Governement
Expendidures, Area 4: PPP deflators, Area 5: Political variables,
Area 6: Trade Policy and others.}}

\newpage

\begin{center}{\footnotesize
\begin{tabular}{rrcc|ccccc}
\hline &M&$\lambda$&$\hat S$&$\widehat{\alpha_2}$&CI
\\\hline
Exp1&BC&1.183&2& -0.016 &[-0.025, -0.004]\\
Exp1&BC&0.788&4& -0.041 &[-0.054, -0.029]\\
Exp1&BC&0.591&3& -0.044 &[-0.065, -0.034]\\
Exp1&BC&0.473&11&-0.051 &[-0.065, -0.032]\\
Exp1&M2&& 3  & -0.136 &[-0.212 , -0.060 ] \\ 
Exp2&M2&&3& -0.116  &[-0.191 , -0.040 ] \\ 
Exp3&M2&&11  & -0.351  &[-0.442 , -0.261 ]\\ 
Exp4&M2&& 11  & -0.332  &[-0.443 , -0.221 ]\\ 
\hline
\end{tabular}

{TABLE 4. Estimation of $\alpha_2$ for selection obtained via lasso
procedure and LOLA for Model2.}}
\end{center}

\newpage

\begin{center}
{\footnotesize
\begin{tabular}{cc|ccc}
\hline &$\hat S$&$\widehat{\alpha_2}$&CI &$N_0$
\\\hline
Exp1& 7.964  & -0.116  &[-0.146 , -0.085 ] & 0.001 \\
Exp2& 7.896  & -0.115  &[-0.144 , -0.085 ] & 0.000 \\
Exp3& 7.958  & -0.130  &[-0.165 , -0.096 ] & 0.002 \\
Exp4& 7.769  & -0.139  &[-0.172 , -0.106 ] & 0.000 \\
&&&&\\
Exp1& 8.772  & -0.085  &[-0.099 , -0.072 ] & 0.002 \\
Exp2& 8.736  & -0.083  &[-0.097 , -0.070 ] & 0.005  \\
Exp3& 8.766  & -0.096  &[-0.111 , -0.080 ] & 0.004 \\
Exp4& 8.773  & -0.094  &[-0.109 , -0.079 ] & 0.002 \\\hline
\end{tabular}}

{\footnotesize TABLE 5. Estimation of $\alpha_2$ for various periods
of time and for predictors $Z$ selected via Model1 (above) and
Model2 (below). All the standard deviation values normalized by
$\sqrt{1000}$ are smaller than $10^{-3}$. }
\end{center}

\newpage

 \begin{center}
 \begin{tabular}{cc}
\includegraphics[width=8cm,height=6cm]{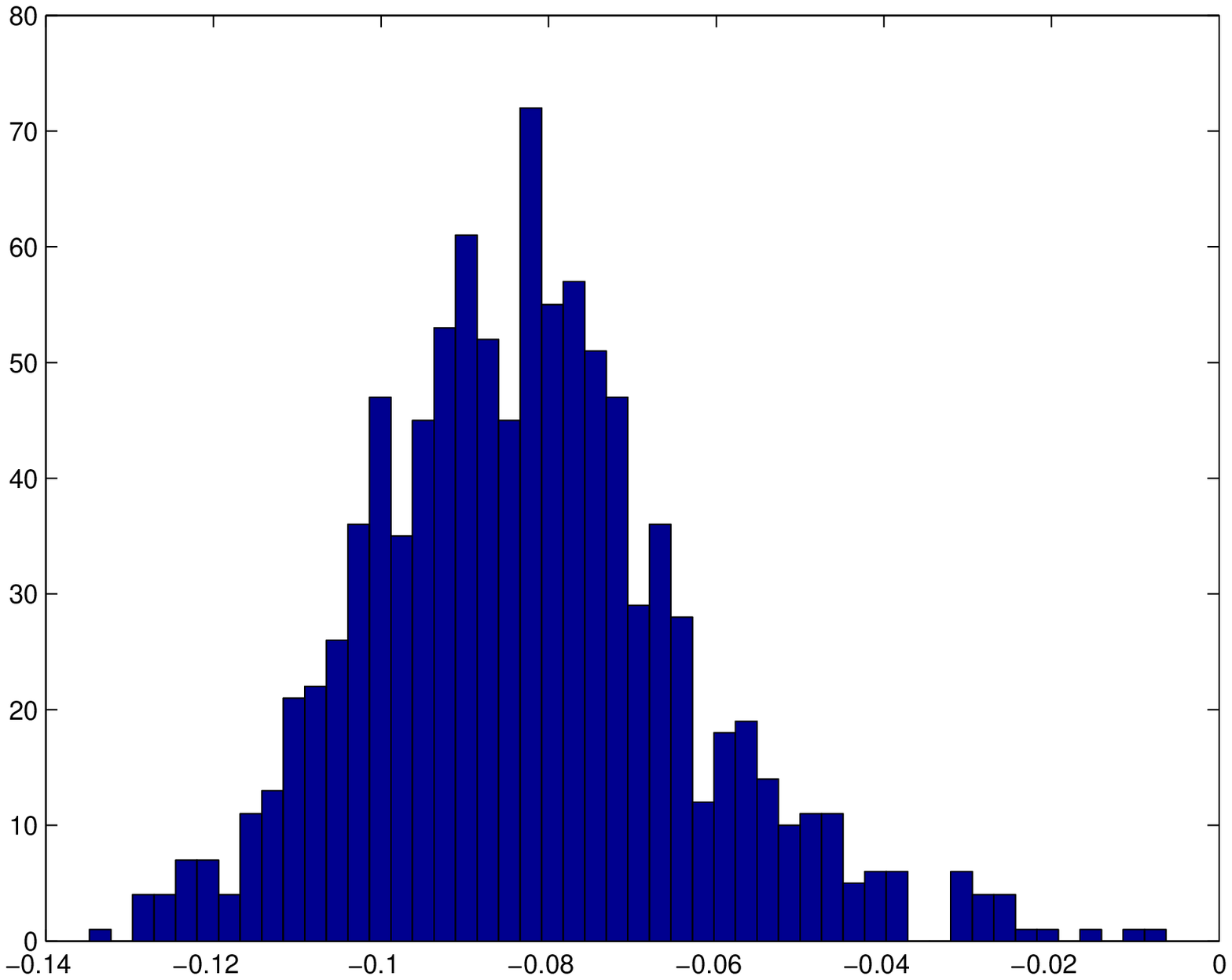} &
\includegraphics[width=8cm,height=6cm]{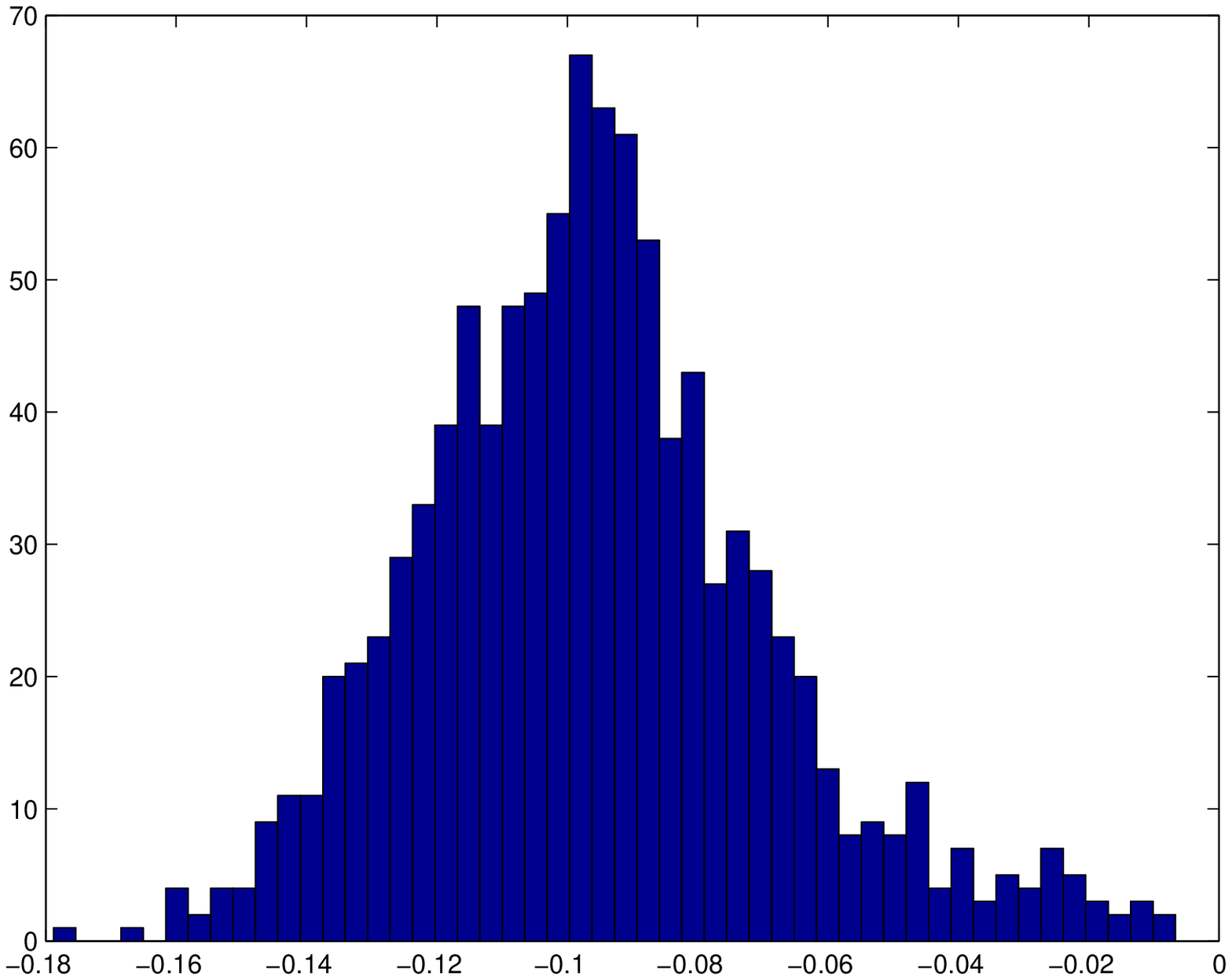} \\
{\footnotesize FIGURE 6. Exp1: 1965-75, $t_0=1960$, empirical
density of $\hat\alpha_2$}&{\footnotesize Exp3: 1975-85, $t_0=1970$,
empirical density of $\hat\alpha_2$}
 \end{tabular}
 \end{center}

\newpage

\begin{center}{\footnotesize
\begin{tabular}{r|ccc|cccc}
\hline &$\hat S$&$\widehat{\alpha_2}$&CI &&$\hat
S$&$\widehat{\alpha_2}$&CI
\\\hline
Exp1& 8  & -0.120  &[-0.148 , -0.092 ] &        &8  & -0.086  &[-0.098 , -0.073 ] \\
Exp2& 8  & -0.116  &[-0.144 , -0.089 ]&         &8  & -0.083  &[-0.095 , -0.071 ] \\
Exp3& 8  & -0.133  &[-0.165 , -0.100 ] &        &9  & -0.100  &[-0.114 , -0.085 ] \\
Exp4&7  & -0.140  &[-0.172 , -0.108 ]  &        &9  & -0.093  &[-0.107 , -0.079 ]  \\
\hline
\end{tabular}

{TABLE 6. Estimation of $\alpha_2$ for both periods of time and for
determinants $Z$ selected via Model1 (left) and Model2 (right). }}
\end{center}

\end{document}